\documentclass[11pt]{article}

\usepackage{amssymb,amsfonts,amsbsy}
\usepackage{amsmath}
\usepackage{psfrag}%
\usepackage{graphicx}
\usepackage{color}
\usepackage{psfrag}
\newtheorem{lemma}{LEMMA}[section]
\newtheorem{theorem}[lemma]{THEOREM}

\newtheorem{corollary}[lemma]{COROLLARY}

\newtheorem{remark}[lemma]{REMARK}

\def\supp{\mathop{\mbox{\rm supp}}\nolimits}
\def\refn#1{{\rm(\ref{#1})}}
\newcommand{\nc}{\newcommand}
\nc{\ha}{\frac{1}{2}}
 \nc{\s}{\widetilde}
  \nc{\dst}{\displaystyle}
\nc{\gm}{\gamma}
 \nc{\ga}{\Gamma}
 \nc{\ka}{\kappa}
\nc{\eps}{\varepsilon}
 \nc{\vep}{\varepsilon}
\nc{\hi}{\varphi} \nc{\vfi}{\varphi}
 \nc{\oa}{\Omega}
\nc{\Om}{\Omega} \nc{\om}{\omega}
 \nc{\ov}{\overline}
 \nc{\lon}{\longrightarrow}
 \nc{\fa}{\forall}
\nc{\scr}{\scriptstyle}
 \nc{\ex}{\exists}
  \nc{\fo}{\forall}
\nc{\pa}{\partial} \nc{\und}{\underline} \nc{\ze}{\zeta}
\nc{\si}{\sigma}
 \nc{\tri}{\triangle}
 \nc{\al}{\alpha}
\nc{\bt}{\beta}
 \nc{\be}{\beta}
  \nc{\ts}{\textstyle}
\nc{\lf}{\left}
\nc{\ra}{\rangle}
\nc{\la}{\langle}
 \nc{\ri}{\right}
 \nc{\lm}{\lambda}
\nc{\lam}{\lambda}
 \nc{\dt}{\delta}
 \nc{\de}{\delta}
\nc{\te}{\theta}
\nc{\vte}{\vartheta}
 \nc{\tl}{\tilde}
 \nc{\wt}{\widetilde}
\nc{\p}{\prime} \nc{\m}{\mu}
 \nc{\R}{{\mathbb{R}}}
 \nc{\B}{{\mathbb{B}}}
 \nc{\N}{{\mathbb{N}}}
 \nc{\C}{{\mathbb{C}}}
 \nc{\bs}{\backslash}

\def\grad{\mathop{\rm grad}\nolimits}
\def\curl{\mathop{\rm curl}\nolimits}
\def\div{\mathop{\rm div}\nolimits}

\def\ker{\mathop{\rm ker}\nolimits}

\def\sup{\mathop{\rm sup}\nolimits}

\def\loc{\mathop{\rm loc}\nolimits}
\def\supp{\mathop{\rm supp}\nolimits}
\def\refn#1{{\rm(\ref{#1})}}
 \nc{\Bc}{{\cal B}}
 \nc{\Ac}{{\cal A}}
 \nc{\Cc}{{\cal C}}
 \nc{\Dc}{{\cal D}}
 \nc{\Ec}{{\cal E}}
 \nc{\Fc}{{\cal F}}
 \nc{\Gc}{{\cal G}}
 \nc{\Hc}{{\cal H}}
 \nc{\Ic}{{\cal I}}
 \nc{\Jc}{{\cal J}}
 \nc{\Kc}{{\cal K}}
 \nc{\Lc}{{\cal L}}
 \nc{\Mc}{{\cal M}}
 \nc{\Nc}{{\cal N}}
 \nc{\Oc}{{\cal O}}
 \nc{\Pc}{{\cal P}}
 \nc{\Qc}{{\cal Q}}
 \nc{\Rc}{{\cal R}}
 \nc{\Sc}{{\cal S}}
 \nc{\Tc}{{\cal T}}
 \nc{\Vc}{{\cal V}}
 \newcommand{\Endbox}{$\hspace*{\fill}\Box$}
 \makeatletter\@addtoreset{equation}{section}\makeatother
  \def\text#1{\;\;\mbox{#1}\;}

 \vspace{
 1cm}
 \makeatletter
 \renewcommand{\@oddhead}{\scriptsize\em Friction.tex
 \hfill    A.Gachechiladze, R.Gachechiladze, J.Gwinner, D.Natroshvili
\hfill \today} \vspace{1cm}

\begin{document}

\vspace{2cm}
\begin{quote}
\begin{center}

 {\bf INTERIOR AND EXTERIOR CONTACT PROBLEMS WITH FRICTION FOR HEMITROPIC
SOLIDS: BOUNDARY VARIATIONAL INEQUALITY APPROACH}
 \end{center}

 \vspace{1cm}

{
  A.~Gachechiladze$^{*}$,  R.~Gachechiladze$^{*}$, J.~Gwinner$^{**}$,
  and D.~Natroshvili \footnote{Corresponding author, { E-mail:
natrosh@hotmail.com}},}

$^{*}${\it A.Razmadze  Mathematical Institute,
 Georgian Academy of Sciences,\\
  M.Aleksidze St. 1, Tbilisi 0193, Republic of Georgia}

$^{**}${\it  Institut f\"ur Mathematik,
Fakult\"at f\"ur Luft- und Raumfahrttechnik,\\ Universit\"at der
Bundeswehr M\"unchen, Werner Heisenberg Weg 39, \\
D-85579 Neubiberg, Germany}

$^{1}${\it Department of Mathematics, Georgian Technical University, \\
M.Kostava St. 77, Tbilisi 0175, Republic of Georgia}

\end{quote}
\date{}

\vspace{1cm}
\begin{quote}
 {\bf Abstract}
We study the interior and exterior contact problems for
hemitropic elastic solids. We treat the cases when the
friction effects, described by
Tresca friction (given friction model), are taken into
consideration either on some part of the boundary of the body or
on the whole boundary. We equivalently reduce these problems to
a boundary variational inequality with the help of the
Steklov-Poincar\'{e} type operator. Based on our boundary
variational inequality approach we prove existence and uniqueness
theorems for weak solutions. We prove that the solutions
continuously depend on the data of the original problem and on the
friction coefficient. For the interior problem necessary and sufficient
conditions of solvability are established when
friction is taken into consideration on the whole boundary.

 \vspace{0.5cm}


{\bf Acknowledgements.} This research was supported by the
Georgian National Science Foundation (GNSF) grant No. GNSF/ST07/3-170.
\end{quote}

\section{Introduction}

{\*}

The main goal of the present paper is the study of contact problems for
hemitropic elastic solids with friction
obeying the Tresca friction model, their mathematical
modelling as
{\it nonsmooth} boundary value problems  and their
analysis with the help of the {\it boundary variational inequality}
technique.

Technological and industrial developments, and also essential
success in biological and medical sciences require to use more
generalized and refined models for elastic bodies.
 In recent years, theories of continuum mechanics  with a complex
 microstructure have been the object of intensive research. Classical
 elasticity associates only the three translational degrees of freedom
 to material points of the body and all the mechanical characteristics
 are expressed by the corresponding displacement vector.
 On the contrary, micropolar theory, by including intrinsic rotations
 of the particles, provides a rather complex model of an elastic body
 that can support body forces and body couple vectors
 as well as  force stress vectors and couple stress vectors at the surface.
  Consequently, in micropolar theory all the mechanical quantities
 are  written in terms of the displacement and microrotation vectors.

 The origin of the rational theories of polar continua
 goes back to brothers E. and F. Cosserat \cite{CC1}, \cite{CC2},
 who gave a development of the mechanics of continuous media
 in which  each material point has the six degrees of freedom defined by 3
 displacement components and 3 microrotation components (for the
 history of  the theory of micropolar elasticity
 see  \cite{Dy1}, \cite{KGBB1}, \cite{Min1}, \cite{Now1}, and the references therein).

A micropolar solid which is not isotropic with respect to
inversion is called {\it hemitropic, noncentrosymmetric,} or {\it
chiral}. Materials may exhibit chirality  on the atomic scale, as
in quartz and in biological molecules - DNA, as well as on a large
scale, as in composites with helical or screw--shaped inclusions,
certain types of nanotubes, bone, fabricated structures such as
foams, chiral sculptured thin films and twisted fibers.
 Experiments have shown that micropolar materials
possess quite different properties in comparison with the
classical elastic materials. For example, the twisting behaviour
under an axial stress is a purely hemitropic (chiral) phenomenon
and  has no counterpart in classical elasticity.  For more details
and applications see the references \cite{AK1}, \cite{AK2}, \cite{CC1},
\cite{Dy1}, \cite{Er1}, \cite{HZ1}, \cite{La1}, \cite{LB1}, \cite{Mu1},
\cite{Mu2}, \cite{Now1}, \cite{Ro1}, \cite{Sh1}, \cite{YL1}.

Refined mathematical models describing the hemitropic  properties
of elastic materials have been proposed by Aero and Kuvshinski
\cite{AK1}, \cite{AK2}. In the mathematical theory of hemitropic
elasticity there are introduced the asymmetric force stress tensor
and moment stress tensor, which are kinematically related with the
asymmetric strain tensor and torsion (curvature) tensor.
The governing equations in
this model become very involved and generate $6\times 6$ matrix
partial differential operator of second order.

 In  \cite{NGGS1}, \cite{NGS1}, \cite{NGZ1}, \cite{NS1} the fundamental matrices
of the associated systems of partial differential equations of
statics and steady state oscillations have been constructed
explicitly in terms of elementary functions and the basic boundary
value and transmission problems of hemitropic elasticity have been
studied by the potential method for smooth and non-smooth
Lipschitz domains. Particular problems of the elasticity theory of
hemitropic continuum  have been considered in \cite{EL1},
\cite{La1}, \cite{LB1}, \cite{LVV1}, \cite{LVV2},  \cite{Now1},
\cite{Now2}, \cite{NN1}, \cite{We1}. The frictionless
unilateral contact problems for hemitropic solids have been studied in \cite{GGN1}.

In classical elasticity similar contact problems have been
considered in many monographs and papers (see,  e.g.,
\cite{DuLi1},
\cite{EJK},
 \cite{Fi1}, \cite{Fi2},  \cite{GaNa1}, \cite{GS1},
\cite{HLNL1},  \cite{HH1},
\cite{Han1} \cite{KiOd1},
\cite{Rod1}, \cite{SST}, and the references therein).

The paper is organized as follows. First we give
the general field equations of the linear theory of elasticity for hemitropic
materials. Then we present a reasonable mathematical model for the
boundary conditions that apply to hemitropic solids in contact
with friction. We start with interior problems and consider the case
when some portion of the
boundary is mechanically fixed and the original problem is
modelled as a coercive boundary variational inequality. Further, we
treat a more complicated case when only traction-contact
conditions are considered on the whole boundary. For this problem
 the corresponding bilinear from is not
coercive any more and we need
the more involved theory of semicoercive variational problems
(see \cite{EJ1}, \cite{Goe}, \cite{GG1} for related problems in classical linear
elasticity and nonlinear elasticity).
In this more involved case,
we establish the necessary and sufficient conditions of solvability.
Next we show that the similar exterior problems are uniquely solvable.
On the basis of the results obtained we prove that
solutions of the boundary variational inequality and, consequently,
the corresponding solutions of the original contact
problems continuously depend on the data of the  problem
and on the friction coefficient.

\setcounter{equation}{0}
\section{Field equations and Green's formulas}
\subsection{Basic Equations}

Let $\Om^{+}\subset\R^3$ be a bounded domain with a $C^\infty$
smooth (we can later relax this assumption), simply connected boundary
$ S:=\pa\Om^{+},\;\ov{\Om^{+}}=\Om^{+}\cup S$. Denote
$\Om^-:=\R^3\setminus\ov{\Om^{+}}$. It is evident that
$\pa\Om^-=S$.

We assume that $\ov{\Om}\in\{\ov{\Om^+},\ov{\Om^-}\}$ is occupied
by a homogeneous hemitropic elastic material. Denote by
$u=(u_1,u_2,u_3)^\top$ and $\om=(\om_1,\om_2,\om_3)^\top$ the {\it
displacement vector} and the {\it micro-rotation vector},
respectively.
Here and in what follows the symbol $(\cdot)^\top$
denotes transposition. Denote by
$n(x)=(n_1(x),n_2(x),n_3(x))^\top$ the
outward normal vector to the surface $S$ at the point $x\in S$.

In hemitropic elasticity theory we have the following
constitutive equations for the {\it force stress tensor}
$\{\tau_{pq} \}$ and the {\it couple stress tensor} $\{\mu_{pq}\}$
for $p,q = 1,2,3$:
\begin{eqnarray}
 &&
 \tau_{pq}=\tau_{pq}(U):=(\mu+\alpha)\, \pa_p u_q +
(\mu-\alpha)\, \pa_q u_p +\lambda\delta_{pq}\div u+
\delta\,\delta_{pq}\div\omega \nonumber \\
 \label{2.1}
 &&
 \hskip2cm +(\varkappa+\nu)\, \pa_p\omega_q+ (\varkappa-\nu)\, \pa_q\omega_p-
2\alpha\sum\limits_{k=1}^{3}\varepsilon_{pqk}\omega_k, \\ 
 &&
 \mu_{pq}=\mu_{pq}(U):=\delta\,\delta_{pq}\div u+
(\varkappa+\nu)\,\Big[ \,\pa_p u_q -
\sum\limits_{k=1}^{3}\varepsilon_{pqk}\omega_k\Big]
+\beta\,\delta_{pq}\div\omega      \nonumber           \\
\label{2.2}
 &&
 \hskip2cm +(\varkappa-\nu) \Big[\,\pa_q u_p -
\sum\limits_{k=1}^{3}\varepsilon_{qpk}\omega_k\Big]
+(\gamma+\varepsilon)\, \pa_p\omega_q +
(\gamma-\varepsilon) \,\pa_q\omega_p , 
\end{eqnarray}
where $U=(u,\om)^\top$,
$\pa=(\pa_1,\,\pa_2,\,\pa_3)$ with $\pa_j=\pa/\pa x_j$,
$\delta_{pq}$ is the Kronecker delta,
$\varepsilon_{pqk}$ is the permutation (Levi-Civit\'a) symbol, and
$\alpha$,  $\beta$, $\gamma$, $\delta$, $\lambda$, $\mu$, $\nu$,
$\varkappa$ and   $\varepsilon$ are the material constants, see
\cite{AK1}, \cite{NGS1}.

The components of the force stress vector
$\tau^{(n)}=(\tau^{(n)}_1,\tau^{(n)}_2,\tau^{(n)}_3)^\top$ and the
couple stress vector
$\mu^{(n)}=(\mu^{(n)}_1,\mu^{(n)}_2,\mu^{(n)}_3)^\top$, acting on
a surface element with the normal vector $n$ read as
\begin{equation}
\label{2.3}
 \tau^{(n)}_q=\sum\limits_{p=1}^3
\tau_{pq}\,n_p,\quad \quad
\mu^{(n)}_q=\sum\limits_{p=1}^3\mu_{pq}\,n_p, \;\;\; q=1,2,3.
\end{equation}
Let us introduce the  $6\times 6$ matrix differential 
stress operator $T(\pa,n)$ \cite{NGS1}
\begin{equation}
              \label{2.4}
T(\pa,n)=\left[
\begin{array}{cc}
T^{(1)}(\pa,n) & T^{(2)}(\pa,n)\\[2mm]
T^{(3)}(\pa,n) & T^{(4)}(\pa,n)
\end{array}\right]_{6\times 6},   \quad\quad
T^{(j)}=\left[ T^{(j)}_{pq} \right]_{3 \times 3}, \;\;\;
j=\ov{1,4},
\end{equation}
\begin{equation}
  \label{2.5}
  \begin{array}{l}
\dst T_{pq}^{(1)}(\pa,n)=(\mu+\alpha)\,\delta_{pq}\,  \pa_n+
(\mu-\alpha)\,n_q\, {\pa_p}+\lambda\,
n_p\, {\pa_q},           \\  
\dst
 T_{pq}^{(2)}(\pa,n)=(\varkappa+\nu)\,\delta_{pq}\, {\pa_n}+
 (\varkappa-\nu)\,n_q\, {\pa_p}+\delta\,n_p\, {\pa_q}-
2\,\alpha\sum\limits_{k=1}^3\varepsilon_{pqk}\,n_k, \\     
\dst T_{pq}^{(3)}(\pa,n)=(\varkappa+\nu)\,\delta_{pq}\, {\pa_n}+
(\varkappa-\nu)\,n_q\, {\pa_p}+\delta\,n_p\, {\pa_q}, \\   
\dst
T_{pq}^{(4)}(\pa,n)=(\gamma+\varepsilon)\,\delta_{pq}\,{\pa_n} +
(\gamma-\varepsilon)\,n_q\, {\pa_p}+\beta\, n_p\,
{\pa_q}-2\,\nu\sum\limits_{k=1}^3\varepsilon_{pqk}\,n_k \,,
\end{array}
\end{equation}
where ${\pa_n}={\pa}/{\pa n}$ denotes the normal derivative.

From the formulas \refn{2.1}, \refn{2.2} and \refn{2.3} it
can be easily checked that
\begin{eqnarray}
\label{2.6}
 \left( \tau^{(n)},\,\mu^{(n)}
\right)^\top=T(\pa,n)\,U.
\end{eqnarray}
The equilibrium equations in the theory of hemitropic
elasticity read as, see \cite{AK1}, \cite{NGS1}
\begin{eqnarray*}
&& \sum\limits_{p=1}^3\pa_p\,\tau_{pq} (x )+ \varrho\, \widetilde{X}^{(1)}_q(x )
=0,\\ 
 && \sum\limits_{p=1}^3\,\pa_p\,\mu_{pq}(x )
+\sum\limits_{l,r=1}^3\,\varepsilon_{qlr}\,\tau_{lr}(x )
+\varrho\,\widetilde{X}^{(2)}_q(x ) =0, \;\;q=1,2,3,
\end{eqnarray*}
where $\varrho$ is the mass density of the elastic material,
while $\widetilde{X}^{(1)}
=(\widetilde{X}^{(1)}_1,\widetilde{X}^{(1)}_2,\widetilde{X}^{(1)}_3)^\top$,
and $\widetilde{X}^{(2)}
=(\widetilde{X}^{(2)}_1,\widetilde{X}^{(2)}_2,\widetilde{X}^{(2)}_3)^\top$  are the
body force and body couple vectors, respectively.

Using the constitutive equations \refn{2.1} and \refn{2.2}
we can rewrite the equilibrium equations in terms of the
displacement and micro-rotation vectors,
\begin{eqnarray}
&& \hskip-15mm (\mu+\alpha)\,\Delta \,u(x
)+(\lambda+\mu-\alpha)\grad\div u(x )+
(\varkappa+\nu)\,\Delta\,\omega(x )\nonumber \\
&&  +(\delta+\varkappa-\nu)\grad\div\omega(x )
+2\alpha\curl\omega(x )+\varrho\, \widetilde{X}^{(1)}(x ) =0, \ \nonumber\\
 &&\hskip-15mm(\varkappa+\nu)\,\Delta
\,u(x )+(\delta+\varkappa-\nu)\grad\div u(x)+ 2\alpha\curl
u(x ) \label{2.7}\\
 &&+(\gamma+\varepsilon)\,\Delta\,\omega(x )
+(\beta+\gamma-\varepsilon)\grad\div\omega(x )+4\nu\curl\omega(x
)\
\nonumber\\ 
&&-4\alpha\,\omega(x )+ \varrho\,\widetilde{X}^{(2)}(x )=0,\nonumber
\end{eqnarray}
where $\Delta=\pa^2_1+\pa^2_2+\pa^2_3$ is the Laplace operator.

Let us introduce the matrix differential operator
given by the left hand side of \refn{2.7}:
\begin{equation}
               \label{2.8}
L(\pa):=\left[
\begin{array}{cc}
L^{(1)}(\pa)  & \;\;L^{(2)}(\pa) \\[3mm]
L^{(3)}(\pa)  &\;\; L^{(4)}(\pa)
\end{array}
\right]_{6\times 6} ,
\end{equation}
\begin{eqnarray*}
\begin{array}{l}
  L^{(1)}(\pa):= (\mu+\alpha)\,\Delta\,I_3+
(\lambda+\mu-\alpha)\,Q(\pa),  \\[2mm]
  L^{(2)}(\pa)=L^{(3)}(\pa):=(\varkappa+\nu)\,\Delta\, I_3+
(\delta+\varkappa-\nu)\,Q(\pa)+2\,\alpha \,R(\pa), \\[2mm]
 L^{(4)}(\pa):=[(\gamma+\varepsilon)\,\Delta
-4\,\alpha]\,I_3+(\beta+\gamma-\varepsilon)\,Q(\pa)+ 4\,\nu
\,R(\pa) \,,
\end{array}
\end{eqnarray*}
where and in the sequel  $I_k$ stands for the $k\times k$ unit
matrix and
$$
Q(\pa):=[\,\pa_k\pa_j\,]_{3\times3},\quad\quad
R(\pa):=\left[\begin{array}{ccc}
0 & -\pa_3 & \pa_2 \\
\pa_3 & 0 & -\pa_1 \\
-\pa_2 & \pa_1 & 0
\end{array}\right]_{3\times 3}\,.
$$
It is easy to see that
$$
 R(\pa)u=\left[\begin{array}{c}
\pa_2u_3  -\pa_3u_2 \\
\pa_3u_1  -\pa_1u_3 \\
\pa_1u_2  -\pa_2u_1
\end{array}\right]=\curl u, \quad  Q(\pa)\,u=\grad
\div u.
$$

Thus \refn{2.7} can be written in matrix form as
$$
L(\pa )U(x)+ X (x)=0  \;\; \text{with}\;\;
  U=(u,\omega)^\top,\;\;\;
X=(X^{(1)},X^{(2)})^\top:=(\varrho\, \widetilde{X}^{(1)},\,\varrho\, \widetilde{X}^{(2)})^\top.
$$
Note that the operator $L(\pa)$ is formally self-adjoint, i.e.,
 $L (\pa)=[L (-\pa)]^\top.$

\subsection{Green's formulas}
For  real-valued vector functions $U=(u,\om)^\top $ and
$U'=(u',\om')^\top$ from the class $[C^2(\ov{\Om^+}\,)]^6$
 the following Green formula holds \cite{NGS1}
 \begin{eqnarray}
\label{2.9}
\int\limits_{\Omega^+}\left[ L(\pa)U \cdot U'
+E(U,U')\right]\,dx= \int\limits_{S} \{T(\pa,n)U\}^+\cdot\{
U'\}^+\,dS,
\end{eqnarray}
where the symbols $\{\,\cdot\,\}^\pm$ denote the one sided limits
(trace operators) on $S$ from $\Om^\pm$ respectively, while
 $E(\cdot\, ,\,\cdot)$ is the bilinear form defined by

\begin{eqnarray}  \label{2.10}
 &&
 \hskip-15mm
 E(U,U')=E(U',U)=\sum\limits_{p,q=1}^3
\{\,(\mu+\alpha)u'_{pq}u_{pq}+(\mu-\alpha)u'_{pq}u_{qp} \nonumber \\
 &&
+(\varkappa+\nu)(u'_{pq}\omega_{pq}+\omega'_{pq}u_{pq})
+(\varkappa-\nu)(u'_{pq}\omega_{qp}+ \omega'_{pq}u_{qp}) +
(\gamma+\varepsilon)\omega'_{pq}\omega_{pq} \nonumber \\[2mm]
 &&
 +(\gamma-\varepsilon)\omega'_{pq}\omega_{qp}
+\delta(u'_{pp}\omega_{qq}+\omega'_{qq}u_{pp})+ \lambda
u'_{pp}u_{qq}+\beta\omega'_{pp}\omega_{qq}\},
\end{eqnarray}
where $u_{pq}$ and $\om_{pq}$  are the so called {\it strain} and
{\it torsion} ({\it curvature}) tensors for hemitropic bodies,
\begin{equation}
\label{2.11} \dst u_{pq}=u_{pq}(U)=\pa_pu_q-
\sum\limits_{k=1}^3\eps_{pqk}\om_k,\;\;\om_{pq}=\om_{pq}(U)=\pa_p\om_q,\;\;p,q=1,2,3.
\end{equation}

Here and in what follows $a\cdot b := a^\top b$ is
the usual scalar product of two vectors $a,b$.
We can generalize Green's formula \refn{2.9} to unbounded
domains.
We say that a vector $U=(u,\om)^\top$ satisfies the decay
condition (Z) at infinity if for sufficiently large $|x|$
$$
 u_j(x)\!=\! \Oc(|x|^{-1}),\;\;\om_j(x)\!=\!\Oc(|x|^{-2}),\;\;
  \frac{\pa u_j(x)}{\pa x_k}\!=\! \Oc(|x|^{-2}), \;\;\frac{\pa
\om_j(x)}{\pa x_k}\!=\! \Oc(|x|^{-2}), \;\;k,j\!=\!1,2,3.
$$
Let $U=(u,\om)^\top $ and $U'=(u',\om')^\top$ belong to the class
$[C^2(\ov{\Om^-})]^6$ and satisfy the decay condition (Z) at infinity. Then the
following Green's formula holds
$$
 \int\limits_{\Omega^-}\left[ L(\pa)U \cdot U'
+E(U,U')\right]\,dx= -\int\limits_{S} \{T(\pa,n)U\}^-\cdot\{
U'\}^-\,dS.
$$
 From formulas \refn{2.10} and \refn{2.11} we get
\begin{eqnarray}
                  \label{2.12}
&& E(U,U')= \frac{3\lambda+2\mu}{3}\left(\div \,u+
\frac{3\delta+2\varkappa}{3\lambda+2\mu}\div \,\omega\right)
\left(\div \,u'+\frac{3\delta+2\varkappa}{3\lambda+2\mu}\div
\,\omega'\right) \nonumber \\
&& \hskip1cm +\frac{1}{3}\left(3\beta+2\gamma-
\frac{(3\delta+2\varkappa)^2}{3\lambda+2\mu}\right)(\div\omega)(\div\omega')
                   \nonumber \\
&& \hskip1cm +\frac{\mu}{2} \sum\limits_{k,j=1,\,k\neq
j}^3\left[\frac{\pa u_k}{\pa x_j} +\frac{\pa u_j}{\pa x_k}+
\frac{\varkappa}{\mu}\left(\frac{\pa\omega_k}{\pa x_j}
+\frac{\pa\omega_j}{\pa x_k}\right)\right]            \nonumber \\[3mm]
&&
 \hskip3cm \times \left[\frac{\pa u'_k}{\pa x_j}
+\frac{\pa u'_j}{\pa x_k}+
\frac{\varkappa}{\mu}\left(\frac{\pa\omega'_k}{\pa x_j}
+\frac{\pa\omega'_j}{\pa x_k}\right)\right]        \nonumber
\end{eqnarray}
\begin{eqnarray}
  &&
  \hskip1cm
 +\,\frac{\mu}{3}\sum_{k,j=1}^3
\left[\frac{\pa u_k}{\pa x_k}-\frac{\pa u_j}{\pa x_j}+
\frac{\varkappa}{\mu}\left(\frac{\pa\omega_k}{\pa x_k}
-\frac{\pa\omega_j}{\pa x_j}\right)\right]  \nonumber\\   [3mm]
&& \hskip3cm
\times \left[\frac{\pa u'_k}{\pa x_k} -\frac{\pa u'_j}{\pa x_j}+
\frac{\varkappa}{\mu}\left(\frac{\pa\omega'_k}{\pa x_k}
-\frac{\pa\omega'_j}{\pa x_j}\right)\right] \nonumber\\   [3mm] &&
 \hskip1cm +\left(\gamma-\frac{\varkappa^2}{\mu}\right)
\sum\limits_{k,j=1,\,k\neq j}^3 \left[\frac{1}{2}
\left(\frac{\pa\omega_k}{\pa x_j}+\frac{\pa\omega_j}{\pa
x_k}\right) \left(\frac{\pa\omega'_k}{\pa
x_j}+\frac{\pa\omega'_j}{\pa x_k}\right)\right. \nonumber \\
[3mm] &&
 \hskip3cm + \left. \frac{1}{3} \left(\frac{\pa\omega_k}{\pa
x_k}-\frac{\pa\omega_j}{\pa x_j}\right)
\left(\frac{\pa\omega'_k}{\pa x_k}-\frac{\pa\omega'_j}{\pa
x_j}\right)\right]  \nonumber \\   [3mm] &&
 \hskip1cm +\alpha\left(\curl\,
u+\frac{\nu}{\alpha}\curl\,\omega -2\,\omega\right) \,\cdot \,
\left(\curl\,
u'+\frac{\nu}{\alpha}\curl\,\omega'-2\,\omega'\right) \nonumber \\
[3mm]
 && \hskip3cm +\left( \varepsilon-\frac{\nu^2}{\alpha}
\right)\,\curl\,\omega \,\cdot \,\curl\,\omega'.
\end{eqnarray}

The necessary and sufficient conditions for the potential energy density
function $E(U,U)$ to be a positive definite quadratic form are the following inequalities
(see \cite{AK2}, \cite{Dy1}, \cite{GGN1})
$$
\begin{array}{l}
\mu>0, \;\;\;\alpha>0,\;\;\; \gamma>0,\;\;\;\varepsilon>0,\;\;\;
\lambda+2\mu>0, \;\;\; \mu\,\gamma-\varkappa^2>0,
\;\;\;\alpha\,\varepsilon-\nu^2>0,
\\[3mm]
(\lambda+\mu) (\beta+\gamma)-
 (\delta+\varkappa)^2>0,\;\;\;
(3\lambda+2\mu)(3\beta+2\gamma)-(3\delta+2\varkappa)^2>0,\\[3mm]
  (\mu+\alpha)(\gamma+\varepsilon)-(\varkappa+\nu)^2>0,\;\;
(\lambda+2\mu)(\beta+2\gamma)-(\delta+2\varkappa)^2>0,\\[3mm]
\mu\,[\, (\lambda+\mu)\,(\beta+\gamma)-(\delta+\varkappa)^2\,]+
(\lambda+\mu)(\mu\,\gamma-\varkappa^2)>0,\\[3mm]
\mu\,[
\,(3\lambda+2\mu)\,(3\beta+2\gamma)-(3\delta+2\varkappa)^2\,]+
(3\lambda+2\mu)\,(\mu\,\gamma-\varkappa^2)>0.
\end{array}
$$
Let us note that, if the condition $3\lambda+2\mu>0$ is fulfilled,
which is very natural in classical elasticity \cite{KGBB1}, 
then the above conditions are equivalent to
the following simultaneous inequalities

$$
\begin{array}{l}
\mu>0, \;\;\;\alpha>0,\;\;\; \gamma>0,\;\;\;\varepsilon>0,\;\;\;
3\lambda+2\mu>0,\;\;\;
 \mu\,\gamma-\varkappa^2>0,
\;\;\;\alpha\,\varepsilon-\nu^2>0,  \\[3mm]
 (\mu+\alpha)(\gamma+\varepsilon)-(\varkappa+\nu)^2>0,\;\;\;
(3\lambda+2\mu)(3\beta+2\gamma)-(3\delta+2\varkappa)^2>0.
\end{array}
$$

For simplicity in what follows we assume that $3\lambda+2\mu>0$
and therefore the above conditions
imply positive definiteness of the energy
density quadratic form $E(U,U)$
with respect to the variables $u_{pq}(U)$ and $\om_{pq}(U)$,
 i.e., there exists  a positive constant $c_0>0$
depending only on the material parameters, such that
\begin{equation}               \label{2.13}
E(U,U) \geq c_0\sum\limits_{p,q=1}^3 \left[u_{pq}^2+\om_{pq}^2\right].
\end{equation}

The following assertion describes the null space of the energy density
quadratic form $E(U,U)$ (see \cite{NGS1}).
\begin{lemma}
Let $U=(u,\om)^\top\in [C^1(\ov{\Om^+})]^6$.
Then $E(U,U)=0$ holds in $\Om^+$, if and only if
$$
u(x)=\chi^{(1)}(x)=[a\times x]+b, \quad \quad
\omega(x)=\chi^{(2)}(x)=a,   \quad x\in \Om^+,
$$
where $a$ and $b$ are arbitrary three-dimensional constant vectors
and the symbol $[\cdot\times\cdot]$ denotes the cross product of
two vectors.
\end{lemma}
Vectors of type $\chi=(\chi^{(1)},\chi^{(2)})^\top=([a\times
x]+b,a)^\top$ are called {\it generalized rigid displacement
vectors}. Observe that a generalized rigid displacement vector
vanishes, i.e. $a=b=0$, if it is zero at a single point.

Throughout the paper $L_p$ $(1\leq p\leq\infty)$,
 $H^s=H^s_2$, $s\in\R,$ denote the
Lebesgue and Bessel potential spaces (see, e.g., \cite{Tr1})
with the norms $\|\,\cdot\,\|_{L_p}$ and
$\|\,\cdot\,\|_{H^s}$, respectively.
Moreover, $L_{p,\,comp}(\Om^-)$ is the subspace
of the $L_p$ functions with compact support in 
the unbounded domain $\Om^-$. 
Denote by ${\cal D}(\Om^\pm)$ the class of
$C^\infty$ functions with support in the domains $\Om^\pm$. If
$M$ is an open proper part of the manifold $S$, i.e., $M\subset
S,\;\,M\neq S$,
 then by $H^s(M)$ we denote the restriction of the space $H^s(S)$ on $M$,
$$
H^s(M):=\{r_{_{M}}\vfi:\;\vfi\in H^s(S)\},
$$
where $r_{_{M}}$ denotes the restriction operator on the set $M$.
Further, let
$$
\wt{H}^s(M):=\{\vfi\in H^s(S):\;\supp\vfi\subset\ov{M}\}.
$$
From the positive definiteness \refn{2.13} of the energy density form $E(\cdot,\cdot)$
it follows that the inequlity

\begin{align}  \label{2.14}
 B(U,U):=\int\limits_{\Om^+ } E(U,U) dx
\geq & \;c_1\; \int\limits_{\Om^+ } \Big\{
\sum\limits_{p,q=1}^3[(\pa_p u_q)^2 +(\pa_p \om_q)^2 ]+
\sum\limits_{q=1}^3[u_q^2+\om_q^2] \Big\}\,dx \nonumber\\
&-c_2 \int\limits_{\Om^+ }
\sum\limits_{q=1}^3[\,u_q^2+\om_q^2\,] \,dx
\end{align}

holds true for an arbitrary real-valued vector function
$U\in[C^1(\ov{\Om^+})]^6$ with some positive constants $c_1,
c_2$ depending only on the material parameters.
By standard limiting arguments we
easily conclude that for any $U\in[H^1(\Om^+)]^6$  there holds the
following Korn's type inequality (cf. \cite{Fi1}, Part I, \S 12,
\cite{Ci1}, \S 6.3)
\begin{equation}
                 \label{2.15}
B(U,U)\geq c_1 \,||U||^2_{[H^1(\Om^+)]^6}
-c_2\,||U||^2_{[L_2(\Om^+)]^6}.
\end{equation}
\begin{remark}
By standard limiting arguments Green's formula \refn{2.9} can be
extended to Lipschitz domains and to vector functions $U,U'\in
[H^1(\Omega^+)]^6$, hence $L(\pa)U \in [L_2(\Omega^+)]^6$ $\mbox{\rm
(see, \rm \cite{Ne1}, \cite{LiMa1})}$,
\begin{equation}
                  \label{2.16}
\int\limits_{\Omega^+ }\left[L(\pa)U\cdot U'+ E(U,U') \right]\,dx=
  \langle \{T(\pa,n)U\}^+
 \, , \,\{U'\}^+
 \rangle_{_{\pa \Omega^+}},
\end{equation}
where  $\langle\, \cdot\, , \, \cdot \, \rangle_{\pa \Omega^+}$
denotes the duality between the spaces $[H^{- 1/2}(\pa
\Omega^+)]^6$ and  $[H^{1/2}(\pa \Omega^+)]^6$, which
extends the inner product in the space $[L_2(\pa \Omega^+)]^6$.
By this relation the generalized trace
 $\{T(\pa,n)U\}^+ \in [H^{- 1/2}(\pa \Omega^+)]^6$
of the stress operator is well-defined.

Analogously, for the unbounded domain $\Om^-$ and for vector
functions $U, U'\in [H^1_{loc}(\Om^-)]^6$, satisfying the decay
condition $(Z)$ along with the imbedding we have
 $$
\int\limits_{\Omega^- }\left[L(\pa)U\cdot U'+ E(U,U')
\right]\,dx=-   \langle \{T(\pa,n)U\}^- \, , \,\{U'\}^-  \rangle_{_{\pa \Omega^-}}.
 $$
\end{remark}
\setcounter{equation}{0}
\section{Contact problems with Tresca friction}
\subsection{Coulomb's law and Tresca friction}
Let the boundary $S$ of the domain $\Om^+$ be divided into two
open, connected and non-overlapping parts $S_1$ and $S_2$ of
positive measure, $S=\ov{S_1}\cup\ov{S_2},$ $S_1\cap S_2=\varnothing$.
Assume that the hemitropic elastic body occupying
the domain $\Om^+$ is in bilateral contact with
a foundation along the subsurface $S_2$, i.e., there is no gap between the body
and the foundation.
Denote by $F(x)=(F_1(x),F_2(x),F_3(x))^\top$
the reaction force stress vector which the foundation exerts on
the hemitropic body at the point $x \in S_2$.

Throughout the paper $F_n$ and $F_s$ stand for the normal and
tangential components of the vector $F$: $F_n=F\cdot n$ and
$F_s=F-(F\cdot n)\,n$. Further, let $\Fc(x)$ be the {\it
friction coefficient} at the point $x$. It is a nonnegative scalar
function which  depends on the geometry of the contacting surfaces
and also on the physical properties of interacting materials.

Coulomb's law describing friction for static bilateral contact reads as
follows (see \cite{DuLi1}):
\textsl{
For the force stress vector $F$ there holds
$$
|F_s(x)| \leq \Fc(x) \,|F_n(x)|\,.
$$
 Moreover, if
$$
|F_s(x)|<\Fc(x)\,|F_n(x)|,
$$
 then the tangential component of the displacement vector
vanishes, $u_s(x)=0,$ \emph{and if}
$$
|F_s(x)|=\Fc(x)\,|F_n(x)|,
$$
 then there exist nonnegative functions  $\lambda_1$
\emph{and} $\lambda_2$  which do not vanish simultaneously
and such that
$$
\lambda_1(x)\,u_s(x)=-\lambda_2(x)\,F_s(x).
$$}

In classical elasticity the fixed point approach to frictional contact as
proposed by Panagiotopoulos \cite{Pan},
 employed in the existence proofs (see \cite{EJK}),
and recently numerically realized in \cite{HKD}
leads to an approximating sequence of contact problems with given friction.
In these approximations, also known as contact problems with Tresca friction
(see \cite{SST}), the unknown normal component is replaced by a
given nonnegative slip stress.

Similarly we replace the unknown normal component $F_n$ by some given force
$F_0$ and replace the above treshold by
\begin{equation}     \label{3.2}
                 g(x):=\Fc(x) \, |F_0(x)| \,.
\end{equation}

\subsection{Pointwise and variational formulation
of the bilateral contact problem}

Consider the equation in the domain $\Omega^+$:
\begin{equation}
 \label{3.1}
 L(\pa)\,U+X=0,
\end{equation}
where $L(\pa)$ is the matrix differential operator given by
the formula \refn{2.8},   $U=(u,\om)^\top$,
$X:=(\rho\,\widetilde{X}^{(1)},\rho\,\widetilde{X}^{(2)})^{\top}
\in [L_2(\Omega^+)]^6$.

A vector-function $U=(u,\om)^\top\in [H^1(\Omega^+)]^6$ is a weak
 solution of equation \refn{3.1} in $\Om^+$ if
$$
 B(U,\Phi)=\int_{\Om^+}X\cdot\Phi\,dx\;\quad\forall\,\Phi\in
[{\cal D}(\Om^+)]^6,
$$
where the bilinear form $B(\cdot\,,\,\cdot)$ is given by formula
\refn{2.14}.

Due to the formulas \refn{2.1}-\refn{2.6} for the force
stress and couple stress vectors we have:
$$
\tau^{(n)}(U)=\Tc U=T^{(1)}u+T^{(2)}\om,\quad \mu^{(n)}(U)=\Mc
U=T^{(3)}u+T^{(4)}\om.
$$
It is clear that
$$
\tau^{(n)}_n(U):=(\Tc U)_n, \quad
\tau^{(n)}_s(U):= (\Tc U)_s=\Tc U-\tau^{(n)}_n(U)\,n.
$$
Further, let $X\in [L_2(\Omega^+)]^6, \;\vfi\in [H^{-1/2}(S_2)]^3,
 \;f\in [H^{1/2}(S_1)]^6,\;F_0\in L_\infty (S_2)$,
and  $\Fc :S_2\rightarrow [0,+\infty)$ be a bounded
measurable function. Thus from  formula  \refn{3.2},
the nonnegative function $g \in L_\infty (S_2)$.

Consider the following bilateral mixed contact problem with friction.

{\it Problem}  (A) ({\it coercive case}). \textsl{Find a weak solution
$U=(u,\om)^\top\in [H^1(\Omega^+)]^6$
 of equation \refn{3.1} satisfying the inclusion
$r_{_{S_2}}\{\tau^{(n)}_s(U)\}^+\in [L_\infty (S_2)]^3$
 and the following conditions:
\begin{eqnarray}
\label{3.3}
  \mbox{\rm (i)}
  &&
  r_{_{S_1}}\{U\}^+=f\quad  \mbox{ on}\quad S_1;\\
  \mbox{\rm (ii)}
&&
  r_{_{S_2}}\{\tau^{(n)}_n(U)\}^+=F_0\quad
 \mbox{  on} \quad S_2; \nonumber\\
\label{3.4}
  \mbox{\rm (iii)}
  &&
    r_{_{S_2}}\{\Mc
U\}^+=\vfi\quad  \mbox{ on}\quad
S_2;\\
\mbox{\rm (iv)}
&&
\hskip-4mm
 \mbox{\rm (a)}\;\; \mbox{\rm if}
\;\;|r_{_{S_2}}\{\tau^{(n)}_s(U)\}^+|<g,\;\; \mbox{ then}\;\; r_{_{S_2}}\{u_s\}^+=0,\nonumber\\
&&
\hskip-4mm
 \mbox{\rm (b)}\;\; \mbox{\rm if} \;\;|r_{_{S_2}}\{\tau^{(n)}_s(U)\}^+|=g, \;\;
\mbox{  then  there  exist nonnegative}\nonumber\\
 &&
 \hskip4mm
 \mbox{functions}\; \lambda_1\;  \mbox{\rm and}\; \lambda_2\;
 \mbox{which do not vanish simultaneously and} \nonumber\\
&&
\hskip4mm \lambda_1\,r_{_{S_2}}\{u_s
\}^+=-\lambda_2\,r_{_{S_2}}\{\tau^{(n)}_s(U)\}^+ \,. \nonumber
 \end{eqnarray}}

We emphasize that by the requirement
$r_{_{S_2}}\{\tau^{(n)}_s(U)\}^+\in [L_\infty (S_2)]^3$ the contact conditions
in (iv) can be understood to hold almost everywhere on $S_2$.
Thus we are here more precise with the pointwise formulation
 than other expositions of contact problems
in classical elasticity (see \cite{KiOd1},\cite{SST}).

 To  reduce Problem (A) to a boundary variational inequality we need
 first to reduce the nonhomogeneous equation \refn{3.1} and
 nonhomogeneous condition \refn{3.3} to the homogeneous
 ones. To this purpose consider the following auxiliary linear mixed
 boundary value problem:

Find a vector function $U_0=(u_0,\om_0)^\top\in
[H^1(\Omega^+)]^6$, which is a weak solution of equation
\refn{3.1} and satisfies the following conditions:
\begin{eqnarray*}
&&
r_{_{S_1}}\{U_0\}^+=f\quad  \mbox{\rm on}\quad S_1;\\
&& r_{_{S_2}}\{T(\pa,n)\,U_0\}^+=0\quad  \mbox{\rm on}\quad S_2.
\end{eqnarray*}
This problem has a unique weak solution \cite{NGS1}. Clearly, if
$W$ is a solution of Problem (A) and $U_0$ is a solution of
the above mixed auxiliary problem, then the difference $U:=W-U_0$
will solve the following problem:

{\it Problem}  ($A_0$).  \textsl{ Find a weak solution $U=(u,\om)^\top\in [H^1(\Omega^+)]^6$ of
the equation
\begin{equation}
 \label{3.5}
 L(\pa)\,U=0\quad \text{in}\quad \Omega^+,
\end{equation}
satisfying the following conditions
$r_{_{S_2}}\{\tau^{(n)}_s(U)\}^+\in [L_\infty (S_2)]^3$ and
\begin{eqnarray}
\label{3.6}
  \mbox{\rm (i)}
  &&
  r_{_{S_1}}\{U\}^+=0\quad  \mbox{ on}\quad S_1;\\
\label{3.7}
  \mbox{\rm (ii)}
  &&
  r_{_{S_2}}\{\tau^{(n)}_n(U)\}^+=F_0\quad
 \mbox{ on} \quad
S_2;\\
\label{3.8}
  \mbox{\rm (iii)}
  &&
   r_{_{S_2}}\{\Mc
U\}^+=\vfi\quad  \mbox{ on}\quad
S_2;\\
\label{3.9}
  \mbox{\rm (iv)}
  &&
  \hskip-4mm
   \mbox{\rm (a)}\;\;
\mbox{if} \;\;|r_{_{S_2}}\{\tau^{(n)}_s(U)\}^+|<g,\;\;
\mbox{ then}\;\; r_{_{S_2}}\{u_s\}^+=\varphi_0,\\
&&
  \hskip-4mm
 \mbox{\rm (b)}\;\;\mbox{if} \;\;|r_{_{S_2}}\{\tau^{(n)}_s(U)\}^+|=g, \;\;
\mbox{ then  there  exist nonnegative}\nonumber\\
 &&
 \hskip4mm
 \mbox{functions}\; \lambda_1\;  \mbox{ and}\; \lambda_2\;
 \mbox{which do not vanish simultaneously and} \nonumber\\
\label{3.10}
&&
\hskip4mm
\lambda_1\,r_{_{S_2}}\{u_s
\}^+=-\lambda_2\,r_{_{S_2}}\{\tau^{(n)}_s(U)\}^++\lambda_1\varphi_0,
 \end{eqnarray}
 where $g$ is defined  by  formula  \refn{3.2} and
 $\varphi_0:=-r_{_{S_2}}\{u_0
\}_s^+\in [ H^{1/2}(S_2)]^3$.}

\subsection{Reduction of Problem ($A_0$) to a boundary
variational inequality} To reduce equivalently Problem ($A_0$) to
a boundary variational inequality we recall that
a solution vector $U=(u,\om)^\top\in [H^1(\Om^+)]^6$ to equation \refn{3.5}
satisfying the Dirichlet boundary condition
$$
 \{U\}^+=h\quad \text{on}\quad S
$$
with $h \in [ H^{1/2}(S)]^6$,  can be uniquely represented
as a single layer potential (see \cite{NGS1})
$$
U(x)=V(\Hc^{-1}\,h)(x):=\int\limits_S\Gamma(x-y)\,(\Hc^{-1}\,h)(y)\,dS_y,
\quad x \in \Om^+ \,,
$$
where $\Gamma$ is a fundamental solution of the operator $L(\pa)$
 and $\Hc$ is the boundary integral operator generated by the trace of the
single layer potential on the boundary $S$ (see the explicit
expression for $\Gamma$ in \cite{GGN1}, \cite{NGS1}):
\begin{equation}
 \label{3.11}
(\Hc h)(x)=\lim\limits_{\Omega^{\pm}\ni z\to x\in S}
  \int\limits_S\Gamma(z-y)\,h(y)\,dS_y=\{V(h)\}^+\equiv\{V(h)\}^-.
\end{equation}
Note that the single layer potential operator $V$ and the integral operator
$\Hc$ have the following mapping properties (see \cite{NGS1})
\begin{eqnarray}
\label{3.12}
 && V  : [ H^{r}(S)]^6\rightarrow
[H^{r+3/2}(\Omega^{+})]^6
\;\;\Big[\;[H^{r}(S)]^6\to [H^{r+3/2}_{\loc}(\Om^-)]^6\;\Big]\;\;\;\fa r\in\R,\\
&&
\nonumber
 \Hc   : [H^{r}(S)]^6\rightarrow [H^{r+1}(S)]^6\quad\fa r\in\R.
 \end{eqnarray}
 Moreover, the operator $\Hc$ is invertible and
 \begin{eqnarray}
\label{3.14}
 \Hc^{-1}: [H^{r}(S)]^6\rightarrow [H^{r-1}(S)]^6\quad\fa r\in\R.
 \end{eqnarray}
 Further, there hold the following limiting relations
 \begin{eqnarray}
 \label{3.15}
\{T(\pa,n)V(h)\}^\pm=(\mp2^{-1}I_6+\Kc)\,h \;\;\text{on} \;\;S,
\end{eqnarray}
where
\begin{eqnarray}
\label{3.16} \Kc\,h(x)=\int\limits_S [T(\pa_x,n(x))\Gamma(x-y)]
h(y)\,dS_y.
\end{eqnarray}
It is shown in \cite{NGS1} that
$$
\mp2^{-1}I_6+\Kc\; :\;
[H^{-1/2}(S)]^6\rightarrow [H^{-1/2}(S)]^6
$$
is a singular integral  operator of normal type with zero index.

 Let $G^+:[ H^{1/2}(S)]^6\rightarrow [H^{1}(\Omega^{+})]^6$
be the operator defined by the formula
\begin{equation}
\label{3.17}
 G^+ h:=V(\Hc^{-1}h).
\end{equation}
It is clear that $L(\pa) G^+ h=0$ in
$\Om^+$ and $\{G^+h\}^+=h$ on $S$.

From the properties of the trace operator and mapping properties
of the  single layer potential operator it follows that there
exist positive numbers $c_1$ and $c_2$, such that
\begin{equation}
\label{3.18}
 c_1\|h\|_{[ H^{1/2}(S)]^6}\leq\|G^+h\|_{[H^{1}(\Omega^{+})]^6}\leq c_2\|h\|_{[
H^{1/2}(S)]^6}\quad\;\; \fa h\in[ H^{1/2}(S)]^6.
\end{equation}
Define the Steklov-Poincar\'{e} type operator
$$
\Ac^+h:=\{T(\pa,n)(G^+h)\}^+=\{T(\pa,n)V(\Hc^{-1}h)\}^+.
$$
The operator $\Ac^+$ is well-defined and due to
\refn{3.17} we have the following representation,
see \refn{3.16} and \refn{3.17}:
\begin{equation}
\label{3.13}
 \Ac^+=\big(-2^{-1}I_6+\Kc\big)\,\Hc^{-1}.
\end{equation}
Denote by $\Lambda(S)$ the set of restrictions on $S$ of rigid
displacement vectors, i.e.,
\begin{equation}
\label{3.19}
 \Lambda(S):= \{\chi(x)=([a\times x]+b,a)^\top,\quad x\in S\},
\end{equation}
where $a$ and $b$ are arbitrary three-dimensional constant
vectors.

With the help of Green's formula \refn{2.16} with
$U=U'=V(\Hc^{-1}h)$, the relations \refn{3.15}, \refn{3.19}
and the uniqueness theorem for the Dirichlet BVP, we infer that
$\ker \Ac^+=\Lambda(S)$.

Now we formulate the following technical lemma describing the
properties of the Steklov-Poincar\'{e} operator.
\begin{lemma}
The following relations are true:
\begin{eqnarray*}
&&
 \mbox{\rm (a)}\;\;\langle\Ac^+h,\eta\rangle_{S}=\langle\Ac^+\eta,h\rangle_{S} \quad
    \fa h \in [H^{1/2}(S)]^6 \text{and} \; \fa\eta\in[H^{1/2}(S)]^6;\\
&&  \mbox{\rm (b)}\;\;\Ac^+:[H^{1/2}(S)]^6\rightarrow
 [H^{-1/2}(S)]^6 \text{ is a continuous operator};\\
 && 
  \mbox{\rm (c)}\;\;\exists \, c\,'>0,c\,''>0\,:\,\langle\Ac^+h,h\rangle_{S}\geq
 c\,'\,\|h\|^2_{[H^{1/2}(S)]^6}- c\,''\|h\|^2_{[L_2(S)]^6}\;\;\quad\fa
 h\in[H^{1/2}(S)]^6;\\
 &&
  \mbox{\rm (d)}\;\;\exists\, c>0\,:\,\langle\Ac^+h,h\rangle_{S}\geq
 c \, \|h\|^2_{[H^{1/2}(S)]^6}\quad\;\;\fa
 h\in[\widetilde{H}^{1/2}(S^*)]^6,\; S^*\subset S;\\
 &&
  \mbox{\rm (e)}\;\;\exists \, c>0\,:\,\langle\Ac^+h,h\rangle_{S}\geq
 c \, \|h-Ph\|^2_{[H^{1/2}(S)]^6}\quad\fa h\in[{H}^{1/2}(S)]^6,
\end{eqnarray*}
where $S^*$ is a proper part of $S$ and $P$ is the operator of the
orthogonal projection $($in the sense of $L_2(S))$ of the space
$[H^{1/2}(S)]^6$ onto the space $\Lambda(S)$.
\end{lemma}
{\it Proof.} (a) Let $h,\eta\in[H^{1/2}(S)]^6$. Taking into
account the equality $L(\pa)\,G^+h=0$, due to Green's formula \refn{2.16} we get the
equality:
$$
\begin{array}{l}
\langle\Ac^+h,\eta\rangle_{_{S}}=\langle\{T(G^+h)\}^+,\eta\rangle_{_{S}}=B(G^+h,G^+\eta)\\[2mm]
=B(G^+\eta,G^+h)=\langle\{T(G^+\eta)\}^+,h\rangle_{S}=\langle\Ac^+\eta,h\rangle_{_{S}}.
\end{array}
$$
The item (b) is evident since
by \refn{3.13}, $\Ac^+$ is the composition of the
continuous operators $\Hc^{-1}$ and $-2^{-1}I_6+\Kc$.
The item (c) can be shown by the following arguments. For arbitrary
$h\in [ H^{1/2}(S)]^6$ with the help of \refn{2.15} we derive
\begin{eqnarray*}
\langle\Ac^+h,h\rangle_{S}=B(G^+h,G^+h)
\geq c_1\, \|G^+h\|^2_{[H^{1}(\Omega^{+})]^6}-
c_2 \, \| G^+h\|^2_{[L_2(\Om^+)]^6}.
\end{eqnarray*}
Since $\{G^+h \}^+=h$ on $S$ by the trace theorem we have
$$
\|h\|_{[ H^{1/2}(S)]^6}\leq c_3\,\|G^+h\|_{[H^{1}(\Omega^{+})]^6},
$$
where $c_3$ is some positive constant independent of $h$.

On the other hand, since $[L_2(S)]^6$ is
continuously embedded into
$[ H^{-1/2}(S)]^6$, then by virtue of the properties \refn{3.12}
and \refn{3.14} we have for $G^+ = V \,\Hc^{-1}$:
$$
\|V(\Hc^{-1}h)\|_{[L_2(\Omega^{+})]^6}\leq c_4
\|\Hc^{-1}h\|_{[H^{-3/2}(S)]^6}\leq c_5
\|h\|_{[H^{-1/2}(S)]^6}\leq c_6 \|h\|_{[L_2(S)]^6}
$$
with some positive constants $c_4,$ $c_5$ and $c_6$ independent of
$h$.
So, finally we obtain that
$$
\langle\Ac^+h,h\rangle_{S}\geq\frac{c_1}{c_3^2}\|h\|^2_{[H^{1/2}(S)]^6}
- c_2 c_6^2
\|h\|^2_{[L_2(S)]^6}.
$$
The item (e) follows from the item (c) (see also Lemma 5.1 in the reference
\cite{GGN1}).
The item (d) follows from the item (e). The lemma is proved. \Endbox

Our goal is to reduce equivalently Problem ($A_0$) to a
boundary variational inequality. To this end, we introduce
 the following convex continuous, but nondifferentiable
functional on the space $[H^{1/2}(S_2)]^3$
\begin{equation}
\label{3.20}
 j(\psi):=\int\limits_{S_2}
g\,|\psi_s-\varphi_0|dS\quad \fa\psi\in [H^{1/2}(S_2)]^3.
\end{equation}
Further, let us define the closed subspace
\begin{equation}
\label{3.21}
{\mathbb{H}}:=\{h=(h^{(1)},h^{(2)})^\top\in
[H^{1/2}(S)]^6\,:\,r_{S_1}h=0\}
\end{equation}
and consider the variational inequality:
  \textsl{Find a vector function $h_0=(h^{(1)}_0,h^{(2)}_0)^\top\in {\mathbb{H}}$ such
 that the following inequality
 \begin{eqnarray}
\label{3.22}
 \langle\Ac^+h_0,h-h_0\rangle_{S}+ j(h^{(1)})-j(h^{(1)}_0)\geq
 \int\limits_{S_2}
F_0\,(h^{(1)}_n-h^{(1)}_{0n})\,dS+\langle\varphi,r_{S_2}(h^{(2)}-h^{(2)}_0)\rangle_{S_2}
\end{eqnarray}
holds for all  $h=(h^{(1)},h^{(2)})^\top\in {\mathbb{H}}$.} \\
Now we show that the variational inequality \refn{3.22} and Problem $(A_0)$
 are equivalent.
\begin{theorem}
The boundary variational inequality \refn{3.22} and Problem
$(A_0)$ are equivalent in the following sense: if
$U\in[H^{1}(\Omega^{+})]^6$ is a solution of Problem $(A_0)$,
then $\{U\}^+\in [H^{1/2}(S)]^6$ is a solution of the variational inequality
\refn{3.22} and vice versa, if $h_0\in {\mathbb{H}}$ is a solution of the variational
inequality \refn{3.22}, then $G^+h_0\in[H^{1}(\Omega^{+})]^6$ is a
weak solution of Problem $(A_0)$.
\end{theorem}
{\it Proof.} Let $U=(u,\om)^\top\in [H^1(\Omega^+)]^6$ be a
solution of Problem $(A_0)$ and
$h_0=(h^{(1)}_0,h^{(2)}_0)^\top:=\{U\}^+$. In accordance with the
definition of the operator $G^+$ we have  $U=G^+h_0$.
We show that the conditions \refn{3.9} and \refn{3.10} yield the
following inequality:
\begin{eqnarray}
\label{3.23}
&
 r_{S_2}\{\tau_s^{(n)}(G^+h_0)\}^+\cdot r_{S_2}(h^{(1)}_s-h^{(1)}_{0s})
 +g(|r_{S_2}h^{(1)}_s-\varphi_0|-|r_{S_2}h^{(1)}_{0s}-\varphi_0|)\geq  0\\[2mm]
 &
\fa h=(h^{(1)},h^{(2)})^\top\in {\mathbb{H}} .\nonumber
\end{eqnarray}
Indeed, let
$$
 \tau_{0s}^{(n)} := \{\tau_s^{(n)}(G^+h_0)\}^+ \;\;\;\text{on}\;\;S.
$$
If $|r_{S_2} \tau_{0s}^{(n)}| < g$, then $r_{S_2}h^{(1)}_{0s} = \varphi_0$,
and the left hand side of \refn{3.23} becomes
\begin{eqnarray*}
 r_{S_2} \tau_{0s}^{(n)} \cdot ( r_{S_2} h^{(1)}_s - \varphi_0 )
+ g |r_{S_2}h^{(1)}_s-\varphi_0| \\[1ex]
 \geq  ( - |r_{S_2} \tau_{0s}^{(n)}| + g ) \,
|r_{S_2}h^{(1)}_s-\varphi_0| \, \geq \, 0 \,.
\end{eqnarray*}
Otherwise, $|r_{S_2} \tau_{0s}^{(n)}| = g$. If in \refn{3.10} $\lambda_1=0$,
then  $\lambda_2 \not= 0$, hence $r_{S_2} \tau_{0s}^{(n)} = 0, g = 0$ and
\refn{3.23} immediately holds. If $\lambda_1 \not= 0$, then
$\lambda := \lambda_2 / \lambda_1 \geq 0$
and the left hand side of \refn{3.23} becomes by \refn{3.10}
\begin{eqnarray*}
&
 r_{S_2} \tau_{0s}^{(n)} \cdot  r_{S_2} h^{(1)}_s
- r_{S_2} \tau_{0s}^{(n)} \cdot ( \varphi_0 - \lambda
r_{S_2} \tau_{0s}^{(n)} ) + g |r_{S_2}h^{(1)}_s-\varphi_0|
- g \lambda | r_{S_2} \tau_{0s}^{(n)}| \\[1ex]
& \geq  - |r_{S_2} \tau_{0s}^{(n)}| |r_{S_2}h^{(1)}_s-\varphi_0|
+ g |r_{S_2}h^{(1)}_s-\varphi_0| \, = \, 0 \,.
\end{eqnarray*}
Integrate the inequality \refn{3.23} over $S_2$ to obtain
$$ \int\limits_{S_2}\{\tau_s^{(n)}(G^+h_0)\}^+\cdot
 (h^{(1)}_s-h^{(1)}_{0s})\,dS+ j(h^{(1)})-j(h^{(1)}_0)\geq 0.
$$
Hence we have
\begin{eqnarray*}
&&
 \int\limits_{S_2}\{\tau_s^{(n)}(G^+h_0)\}^+\cdot
 (h^{(1)}_s\!-\!h^{(1)}_{0s})\,dS\!+\! j(h^{(1)})\!-\!j(h^{(1)}_0)\!+\!\int\limits_{S_2}\{\tau_n^{(n)}(G^+h_0)\}^+
 (h^{(1)}_n\!-\!h^{(1)}_{0n})\,dS\\
 &&
 +\langle r_{S_2}\{\Mc(G^+h_0)\}^+, r_{S_2}(h^{(2)}\!-\!h^{(2)}_0)\rangle_{S_2}\geq\int\limits_{S_2}
F_0\,(h^{(1)}_n\!-\!h^{(1)}_{0n})\,dS\!+\!\langle\varphi,r_{S_2}(h^{(2)}\!-\!h^{(2)}_0)\rangle_{S_2}.
\end{eqnarray*}
Since $r_{S_1}(h-h_0)=0$ due to the inclusion $h,h_0\in {\mathbb{H}}$, then
finally we arrive at the inequality
\begin{eqnarray*}
 &
 \dst
 \langle\Ac^+h_0,h-h_0\rangle_{S}+ j(h^{(1)})-j(h^{(1)}_0)\geq\int\limits_{S_2}
F_0\,(h^{(1)}_n-h^{(1)}_{0n})\,dS+\langle\varphi,r_{S_2}(h^{(2)}-h^{(2)}_0)\rangle_{S_2} \\
&
 \fa h=(h^{(1)},h^{(2)})^\top\in {\mathbb{H}} .
\end{eqnarray*}
Thus, $\{U\}^+$ solves the variational inequality \refn{3.22}.

Now, let $h_0=(h^{(1)}_0,h^{(2)}_0)^\top\in {\mathbb{H}} $ be a solution of
the variational inequality \refn{3.22}. Denote $U:=G^+h_0$ in
$\Om^+$. We have to show that $U$ solves Problem $(A_0)$. Due
to the definition of the operator $G^+$ the vector $G^+h_0$ is a
weak solution of the equation \refn{3.5} and
$r_{S_1}\{U\}^+=r_{S_1}h_0 =0$, since $h_0\in {\mathbb{H}} $. Thus \refn{3.5}
and \refn{3.6} hold.

Let $h=(h^{(1)},h^{(2)})^\top\in {\mathbb{H}} $ and
$h^{(1)}_{s}=h^{(1)}_{0s},\;h^{(2)}=h^{(2)}_0$ and
$h^{(1)}_n=h^{(1)}_{0n}\pm\psi$, with arbitrary
$\psi\in\wt{H}^{1/2}(S_2)$. Then $j(h^{(1)})=j(h^{(1)}_0)$ and
from \refn{3.22} we obtain
$$
\langle
r_{S_2}\{\tau_n^{(n)}(G^+h_0)\}^+,r_{S_2}\psi\rangle_{S_2}=\int\limits_{S_2}
F_0\,\psi\,dS\quad\fa\psi\in\wt{H}^{1/2}(S_2),
$$
implying
\begin{eqnarray}
\label{3.25}
 r_{S_2}\{\tau_n^{(n)}(G^+h_0)\}^+=F_0\,\quad \text{on}\quad S_2.
\end{eqnarray}
Thus \refn{3.7} holds.

Select $h=(h^{(1)},h^{(2)})^\top\in {\mathbb{H}} $ such that
$h^{(1)}=h^{(1)}_0$ and $h^{(2)}=h^{(2)}_0\pm\psi$, with arbitrary
$\psi\in[\wt{H}^{1/2}(S_2)]^3$. From \refn{3.22} we have
$$
\langle r_{S_2}\{\Mc(G^+h_0)\}^+,r_{S_2}\psi\rangle_{S_2}=\langle
\varphi,r_{S_2}\psi\rangle_{S_2}\quad\fa\psi\in[\wt{H}^{1/2}(S_2)]^3 \,,
$$
hence
\begin{eqnarray}
\label{3.26}
 r_{S_2}\{\Mc U\}^+= r_{S_2}\{\Mc(G^+h_0)\}^+=\varphi\quad \text{on}\quad
 S_2,
\end{eqnarray}
and \refn{3.8} holds.\\
From \refn{3.22}, by virtue of \refn{3.25} and \refn{3.26} we
derive
\begin{eqnarray*}
&
 \dst
 \langle
r_{S_2}\{\tau_s^{(n)}(G^+h_0)\}^+,r_{S_2}(h^{(1)}_{s}-h^{(1)}_{0s})\rangle_{S_2}+\int\limits_{S_2}
g\,(|h^{(1)}_s-\varphi_0|-|h^{(1)}_{0s}-\varphi_{0}|)\,dS\geq 0 \\
 &
 \fa h=(h^{(1)},h^{(2)})^\top\in {\mathbb{H}} .
\end{eqnarray*}
We rewrite this inequality as follows
\begin{eqnarray}
\label{3.27}
&&
 \langle
r_{S_2}\{\tau_s^{(n)}(G^+h_0)\}^+,r_{S_2}h^{(1)}_s-\varphi_0\rangle_{S_2}
- \langle r_{S_2}\{\tau_s^{(n)}(G^+h_0)\}^+,r_{S_2}h^{(1)}_{0s}-\varphi_0\rangle_{S_2}
\nonumber\\[2mm]
&&
 +\int\limits_{S_2} g\,|h^{(1)}_s-\varphi_0|\,dS-\int\limits_{S_2} g\,|h^{(1)}_{0s}-\varphi_0|\,dS\geq 0
\quad\fa h=\big (h^{(1)},h^{(2)}\big )^\top\in {\mathbb{H}} .
\end{eqnarray}
We set
$$
\theta_s:=r_{S_2}h^{(1)}_{s}-\vfi_0,\quad \theta_{0s}:=r_{S_2}h^{(1)}_{0s}-\vfi_0.
$$
Then \refn{3.27} reads
\begin{eqnarray}
\label{3.28}
&&
 \langle
r_{S_2}\{\tau_s^{(n)}(G^+h_0)\}^+,\theta_s\rangle_{S_2}+
\int\limits_{S_2} g\,|\theta_s|\,dS- \langle
r_{S_2}\{\tau_s^{(n)}(G^+h_0)\}^+,\theta_{0s}\rangle_{S_2} \nonumber\\
&&
 -\int\limits_{S_2} g\,|\theta_{0s}|\,dS\geq 0.
\end{eqnarray}
For arbitrary $\psi\in[\wt{H}^{1/2}(S_2)]^3$ we have
$$
\langle
r_{S_2}\{\tau_s^{(n)}(G^+h_0)\}^+,r_{S_2}\psi_s\rangle_{S_2}=\langle
r_{S_2}\{\tau_s^{(n)}(G^+h_0)\}^+,r_{S_2}\psi\rangle_{S_2}
$$
and $|r_{S_2}\psi_s|\leq |r_{S_2}\psi|$. Therefore, if we take
$r_{S_2}\psi_s$ for $\theta_s$ in \refn{3.28}, we obtain
\begin{eqnarray}
\label{3.29}
&&
 \langle
r_{S_2}\{\tau_s^{(n)}(G^+h_0)\}^+,r_{S_2}\psi\rangle_{S_2}+
\int\limits_{S_2} g\,|\psi|\,dS-\{\langle
r_{S_2}\{\tau_s^{(n)}(G^+h_0)\}^+,\theta_{0s}\rangle_{S_2} \nonumber\\
&&
 +\int\limits_{S_2} g\,|\theta_{0s}|\,dS\}\geq 0 \quad \fa\psi\in[\wt{H}^{1/2}(S_2)]^3.
\end{eqnarray}
Put $\pm t\,\psi$ with $t\geq 0$ for $\psi$ in \refn{3.29}:
\begin{eqnarray*}
&& t\,\Big \{\pm \langle
r_{S_2}\{\tau_s^{(n)}(G^+h_0)\}^+,r_{S_2}\psi\rangle_{S_2}+
\int\limits_{S_2} g\,|\psi|\,dS\Big\}-\{\langle
r_{S_2}\{\tau_s^{(n)}(G^+h_0)\}^+,\theta_{0s}\rangle_{S_2} \nonumber\\
&&
 +\int\limits_{S_2} g\,|\theta_{0s}|\,dS\}\geq 0 \quad \fa t\geq 0,
\;\fa\psi\in[\wt{H}^{1/2}(S_2)]^3.
\end{eqnarray*}
First, sending $t$ to $+\infty$ and afterwards sending $t$ to $0$,
we arrive at the inequalities:
\begin{eqnarray}
&&
\label{3.30}
 |\langle
r_{S_2}\{\tau_s^{(n)}(G^+h_0)\}^+,r_{S_2}\psi\rangle_{S_2}|\leq
\int\limits_{S_2} g\,|\psi|\,dS \quad \fa\psi\in[\wt{H}^{1/2}(S_2)]^3, \\
&&
\label{3.31}
\int\limits_{S_2} \big( \{\tau_s^{(n)}(G^+h_0)\}^+ \cdot \theta_{0s}
+ g\,|\theta_{0s}|\big)\,dS\leq 0.
\end{eqnarray}
Consider the linear functional $\Phi$ on the space
$[\wt{H}^{1/2}(S_2)]^3$ given by
$$
\Phi(\psi)=\langle
r_{S_2}\{\tau_s^{(n)}(G^+h_0)\}^+,r_{S_2}\psi\rangle_{S_2}\,.
$$
Due to the inequality \refn{3.30} this functional is continuous on
the space $[\wt{H}^{1/2}(S_2)]^3$ with respect to the topology
induced by the space $[L_1(S_2)]^3$. Since the space
$[\wt{H}^{1/2}(S_2)]^3$ is dense in $[L_1(S_2)]^3$, so the
functional $\Phi$ can be continuously extended to the whole space
$[L_1(S_2)]^3$ with the same norm.
Since the dual of $[L_1(S_2)]^3$ is isomorphic to $[L_\infty(S_2)]^3$,
there exists a function $\Phi^*\in[L_\infty(S_2)]^3$ such that
$$
\Phi(\psi)=\int\limits_{S_2}\Phi^*\cdot\psi\,dS \quad
\fa\psi\in[L_1(S_2)]^3.
 $$
Hence

$$
r_{S_2}\{\tau_s^{(n)}(G^+h_0)\}^+=\Phi^*\in[L_\infty(S_2)]^3.
$$
Using again the inequality \refn{3.30} we derive
\begin{eqnarray}
\label{3.32}
 \int\limits_{S_2}\big[
\pm\{\tau_s^{(n)}(G^+h_0)\}^+\cdot\psi-g\,|\psi|\,\big]\,dS\leq
0\quad\quad \fa\psi\in[\wt{H}^{1/2}(S_2)]^3
\end{eqnarray}
and the
inequality
\begin{eqnarray}
\label{3.33}
 |r_{S_2}\{\tau_s^{(n)}(G^+h_0)\}^+|\leq g
\quad\text{ almost everwhere on} S_2
\end{eqnarray}
follows.

Indeed, it is well known that for an arbitrary essentially bounded
function $\widetilde{\psi}\in L_{\infty}(S_2)$ there is a sequence
$\widetilde{\varphi}_l \in C^{\infty}(S_2)$ with
$\supp\,\widetilde{\varphi}_l\subset S_2$, such that (see, e.g.,
\cite{Ni1}, Lemma 1.4.2)
$$
\begin{array}{l}
\lim\limits_{l \to
\infty}\widetilde{\varphi}_l(x)=\widetilde{\psi}(x)
\text{for almost all} x\in S_2 \text{and} \\[2mm]
|\widetilde{\varphi}_l(x)|\leq \mbox{\rm ess}\sup\limits_{y\in
S_2} |\widetilde{\psi}(y)| \text{for almost all} x\in S_2\,.
\end{array}
$$
Therefore from inequality \refn{3.32} by the Lebesgue dominated
convergence theorem  it follows that
$$
\int\limits_{S_2}\big[\, \pm \{\tau^{(n)}_s(G^+h_0)\}^+ \cdot \psi -
g\,|\psi|\,\big]\, dS \leq 0 \qquad \forall \psi\in
[L_\infty(S_2)]^3.
$$
In the place of $\psi$ we can put here $\chi(S_2^*)\,\psi$ where
$\psi\in [L_\infty(S_2)]^3$ and $\chi(S_2^*)$ is the
characteristic function of an arbitrary measurable subset
$S_2^*\subset S_2$. As a result we arrive at the inequality $\pm
\{\tau^{(n)}_s(G^+h_0)\}^+ \cdot \psi - g\,|\psi|   \leq 0$
almost everywhere on $S_2$ for
all $\psi\in [L_\infty(S_2)]^3$ and consequently by choosing
$\psi=\{\tau^{(n)}_s(G^+h_0)\}^+$  we finally get \refn{3.33}.

By virtue of \refn{3.31} and \refn{3.33} we obtain
\begin{eqnarray*}
& \dst\int\limits_{S_2} g\,|\theta_{0s}|\,dS \leq
- \int\limits_{S_2}  \{\tau_s^{(n)}(G^+h_0)\}^+ \cdot \theta_{0s}\,dS
  \\[1ex]
& \dst\leq \int\limits_{S_2} |\{\tau_s^{(n)}(G^+h_0)\}^+| \,|\theta_{0s}|\,dS
\leq  \int\limits_{S_2} g\,|\theta_{0s}|\,dS \,,
\end{eqnarray*}
hence \refn{3.31} holds with the equality sign, further by \refn{3.33},
the integrand in \refn{3.31} is nonnegative. Thus we arrive at


\begin{eqnarray}
\label{3.34}
 r_{S_2}\{\tau_s^{(n)}(G^+h_0)\}^+\cdot
 \theta_{0s}+g|\theta_{0s}|=0.
\end{eqnarray}
If $|r_{S_2}\{\tau_s^{(n)}(G^+h_0)\}^+|<g$, then  \refn{3.34} yields
$\theta_{0s}=0$. Hence $r_{S_2}h^{(1)}_{0s}=\vfi_0$ that is  $r_{S_2}\{u_s\}^+=\vfi_0$ and
\refn{3.9} holds.
Otherwise, $|r_{S_2}\{\tau_s^{(n)}(G^+h_0)\}^+|=g$. Then we can rewrite
\refn{3.34} as
$$
g\,|\theta_{0s}|\,(\cos\al +1)=0,
$$
where $\al$ is the angle between the vectors
$r_{S_2}\{\tau_s^{(n)}(G^+h_0)\}^+(x)$ and
$r_{S_2}\{\vartheta_{0s}\}^+(x)$, $x\in S_2$. Now, it is clear that
there are functions $\lambda_1(x)\geq 0$ and $\lambda_2(x)\;\geq
0$ with $\lambda_1(x)+\lambda_2(x)>0$, such that
$$
\lambda_1\,\theta_{0s}=-\lambda_2\,r_{S_2}\{\tau_s^{(n)}(G^+h_0)\}^+,
$$
i.e.,
$$
\lambda_1\, r_{S_2}h^{(1)}_{0s}=-\lambda_2\,r_{S_2}\{\tau_s^{(n)}(G^+h_0)\}^++\lambda_1\vfi_0.
$$
Therefore
$$
\lambda_1\,
r_{S_2}\{u_s\}^+=-\lambda_2\,r_{S_2}\{\tau_s^{(n)}(G^+h_0)\}^++\lambda_1\,\vfi_0,
$$
and \refn{3.10} holds. The proof is complete. \Endbox

\section{Existence and  uniqueness of solutions and their dependence on the
data of the problem}

\subsection{Uniqueness}
We start with the following uniqueness result.
\begin{theorem}
Suppose, the Dirichlet boundary part $S_1$ has positive measure. Then the
boundary variational inequality \refn{3.22} has at most one
solution.
\end{theorem}
{\it Proof.} Let $h_0=(h^{(1)}_0,h^{(2)}_0)^\top\in {\mathbb{H}} $ and
$h^*_0=(h^{(1)*}_0,h^{(2)*}_0)^\top\in {\mathbb{H}} $ be two arbitrary
solutions of the variational inequality \refn{3.22}. Then from
\refn{3.22} we have:
\begin{eqnarray*}
&&
 \langle\Ac^+h_0,h^*_0-h_0\rangle_{S}+ j(h^{(1)*}_0)-j(h^{(1)}_0)\geq\int\limits_{S_2}
F_0\,(h^{(1)*}_{0n}-h^{(1)}_{0n})\,dS+\langle\varphi,r_{S_2}(h^{(2)*}_0-h^{(2)}_0)\rangle_{S_2},\\
&&
 \langle\Ac^+h^*_0,h_0-h_0^*\rangle_{S}+ j(h^{(1)}_0)-j(h^{(1)*}_0)\geq\int\limits_{S_2}
F_0\,(h^{(1)}_{0n}-h^{(1)*}_{0n})\,dS+\langle\varphi,r_{S_2}(h^{(2)}_0-h^{(2)*}_0)\rangle_{S_2}.
\end{eqnarray*}
By summing up these inequalities we obtain
$$
 \langle\Ac^+(h_0-h^*_0),h_0-h_0^*\rangle_{S}\leq 0 \,.
$$
Hence, in view of the non-negativity of the operator $\Ac^+$
$$
 \hskip+0,6cm \langle\Ac^+(h_0-h^*_0),h_0-h_0^*\rangle_{S}= 0\,.
$$
By \refn{2.16} we get
\begin{eqnarray*}
&&
0=\langle\{T(\pa,n)V(\Hc^{-1}(h_0-h_0^*))\}^+,h_0-h_0^*\rangle_S=\la
\{T(\pa,n)G^+(h_0-h_0^*)\}^+,\{G^+(h_0-h_0^*)\}^+\ra_S\\[2mm]
&&
\;\;\;=B(V(\Hc^{-1}(h_0-h_0^*)),V(\Hc^{-1}(h_0-h_0^*))).
\end{eqnarray*}
Thus
$$
G^+(h_0-h_0^*)=V(\Hc^{-1}(h_0-h_0^*))=([a\times
x]+b,a)\quad\text{in}\quad\Om^+.
$$
Since $h_0,\;h_0^*\in \mathbb{H} $ we have
$$
r_{S_1}\{V\big(\Hc^{-1}(h_0-h_0^*)\big)\}=r_{S_1}(h_0-h_0^*)=0,\;\;\text{i.e.,}\;\;([a\times
x]+b,a)=0\;\;\text{on}\;\;S_1.
$$
Therefore $a=b=0$ and $V(\Hc^{-1}(h_0-h_0^*))=0$ in $\Om^+.$
Hence we conclude that $h_0=h_0^*$ on $S$. The proof is complete.
\Endbox

\subsection{Existence results}
Consider the following functional on the closed subspace ${\mathbb{H}}$
(see \refn{3.21})
\begin{eqnarray} \label{4.1}
 \hskip-5mm
 {\cal J}(h)=\frac{1}{2}\langle\Ac^+h,h\rangle_{S}+ j(h^{(1)})-\int\limits_{S_2}
F_0\,h^{(1)}_n\,dS+\langle\varphi,r_{S_2}h^{(2)}\rangle_{S_2} \quad
\fa h=(h^{(1)},h^{(2)})^\top\in {\mathbb{H}}.
\end{eqnarray}
It is easy to show that, due to the self-adjointness property of
the operator $\Ac^+$ (see Lemma\,3.1.\,(a)), the solvability of
the boundary variational inequality \refn{3.22} is equivalent to
the minimization problem for the functional \refn{4.1} on the set
${\mathbb{H}}$.

Since $j(h^{(1)})\geq 0$ and  the operator $\Ac^+$ is bounded  from below
 on ${\mathbb{H}}$ (see Lemma\,3.1.\,(d)) we have
$$
{\cal J}(h)\geq
c_1\|h\|^2_{[H^{1/2}(S)]^6}-c_2\|h\|_{[H^{1/2}(S)]^6}\quad \fa
h\in {\mathbb{H}} .
$$
Consequently, when $h\in {\mathbb{H}}$ and $\|h\|_{[H^{1/2}(S)]^6}\rightarrow
\infty$, then ${\cal J}(h)\rightarrow +\infty$. Therefore the functional
${\cal J}$ in \refn{4.1} is coercive on the closed subspace ${\mathbb{H}}$. Moreover, the
functional ${\cal J}$ is convex and continuous. Due to the general theory
of variational inequalities (see \cite{DuLi1},
\cite{GLT1}) we conclude that the variational inequality \refn{3.22}  has a unique solution.
  In turn this
implies the existence and uniqueness theorems for  Problems $(A_0)$ and  $(A)$.
\begin{theorem}
Suppose, the Dirichlet boundary part $S_1$ has positive measure.
Let $\vfi\in [H^{-{1/2}}(S_2)]^3,\;F_0\in
L_\infty(S_2),\;\vfi_0\in[H^{1/2}(S_2)]^3$. Then Problem
$(A_0)$ is uniquely solvable in the space $[H^1(\Om^+)]^6$
and the solution is representable in the form $U=G^+h_0$, where $h_0$ is a unique solution of the
variational inequality \refn{3.22}.
\end{theorem}
{\it Proof.} It immediately follows from Theorem 3.2
and Theorem 4.1.  \Endbox

\begin{corollary} Let $X\in[L_2(\Omega^+)]^6,$ $\vfi\in
[H^{-{1/2}}(S_2)]^3,$ $F_0\in L_\infty(S_2),$ $f\in
[H^{1/2}(S_1)]^6$ and $\Fc:\;S_2\rightarrow[0,\infty)$ be a
bounded measurable function. Then Problem $(A)$ has a unique
solution in the space $[H^1(\Om^+)]^6$.
\end{corollary}

\subsection{ Lipschitz continuous dependence
of solutions on the problem data}
Let $U\in[H^1(\Om^+)]^6$ and $\widetilde{U}\in[H^1(\Om^+)]^6$ be two
solutions of Problem $(A_0)$ corresponding to the data
$F_0,\;\vfi,\;g$ and $ \widetilde{F}_0,\; \widetilde{\vfi},\; \widetilde{g}$
respectively. Further, let
$h_0=(h^{(1)}_0,h^{(2)}_0)^\top\in {\mathbb{H}} \subset [H^{1/2}(S)]^6$ and
${\widetilde{h}}_0=( {\widetilde{h}}^{(1)}_0,  {\widetilde{h}}^{(2)}_0)^\top\in {\mathbb{H}} \subset [H^{1/2}(S)]^6$ be
the traces of the vector-functions $U$ and $\widetilde{U}$ on the surface
$S$. Then, by virtue of Theorem 3.2, the vectors $h_0$ and ${\widetilde{h}}_0$
will be two solutions of the variational inequality \refn{3.22}
corresponding to the above data. So we have two variational
inequalities of the type \refn{3.22}, the first one for ${h}_0$ and
the second one for ${\widetilde{h}}_0$ . Substitute $h={\widetilde{h}}_0$ in the first
one and $h=h_0$ in the second one, and
sum up to obtain
\begin{eqnarray*}
&
\dst
 -\langle\Ac^+(h_0- {\widetilde{h}}_0),h_0-  {\widetilde{h}}_0\rangle_{S}
 -\int\limits_{S_2}(g -  {\widetilde{g}})(|h_{0s}^{(1)}
 -\vfi_0|-|  {\widetilde{h}}_{0s}^{(1)}-\vfi_0|)\,dS\\
&
\dst
 \geq -\int\limits_{S_2}(F_0-  {\widetilde{F}_0})(h_{0n}^{(1)}-{\widetilde{h}}_{0n}^{(1)})\,dS
 - \langle\vfi-{\widetilde{\vfi}},r_{S_2}({h}^{(2)}_{0s}- {\widetilde{h}}^{(2)}_{0s})\rangle_{S}.
\end{eqnarray*}
From this inequality, taking into account \refn{3.18} and also the
property (d) of the operator $\Ac^+$, we can easily derive the
following Lipschitz estimate:
\begin{eqnarray*}
&
 \|U- {\widetilde{U}}\|_{[H^1(\Om^+)]^6}\leq c_1\,\|h_0- {\widetilde{h}}_0\|_{[H^{1/2}(S)]^6} \\
 &
 \leq c_2\,\Big(\|g- {\widetilde{g}}\|_{L_2(S_2)}+\|F_0 -  {\widetilde{F}_0}\|_{L_2(S_2)}+
\|\vfi-  {\widetilde{\vfi}}\|_{[H^{-1/2}(S)]^3}\Big),
\end{eqnarray*}
where the positive constants $c_1$ and $c_2$ do not depend on the data of the problem.

\section{The semicoercive case}

Let $S_1=\varnothing$. Then $S_2=S$ and for the corresponding
Problem (A) we will have the following formulation.
 Assume that $X\in[L_2(\Omega^+)]^6,$ $F_0\in
L_\infty(S),$ $\vfi\in [H^{-{1/2}}(S)]^3,$ $\Fc:\;S\rightarrow[0,\infty)$ is a bounded
measurable function and $g=\Fc|F_0|$.

{\it Problem $(B)$ $($semicoercive case$)$}.
 \textsl{Find a vector-function
$U=(u,\om)^\top\in [H^1(\Omega^+)]^6$ which is a weak solution of the equation
\begin{equation}
 \label{5.1}
 L(\pa)\,U+X=0\; \text{in}\; \Om^+,
\end{equation}
 satisfying  on $S$ the  conditions $\{\tau^{(n)}_s(U)\}^+\in [L_\infty (S)]^3$ and
\begin{eqnarray}
 \mbox{\rm (i)}
 &&
  \{\tau^{(n)}_n(U)\}^+=F_0; \nonumber\\
 \mbox{\rm (ii)}
 &&
 \{\Mc U\}^+=\vfi;\nonumber\\
 \mbox{\rm (iii)}
 &&
 \hskip-4mm
  \mbox{\rm(a)}\;\; \mbox{if}\; \;|\{\tau^{(n)}_s(U)\}^+|<g,\;\; \mbox{then}\;\; \{u_s\}^+=0, \nonumber\\
&&
\hskip-4mm
 \mbox{\rm (b)}\;\; \mbox{if}\; \;|\{\tau^{(n)}_s(U)\}^+|=g, \;\; \mbox{then there exist nonnegative}\nonumber\\
 &&
 \hskip4mm
 \mbox{functions}\; \lambda_1\;  \mbox{and}\; \lambda_2\;
 \mbox{which do not vanish simultaneously and} \nonumber\\
 &&\hskip4mm \lambda_1\,\{u_s
\}^+=-\lambda_2\,\{\tau^{(n)}_s(U)\}^+.\nonumber
 \end{eqnarray}}

Let the boundary $S$ of $\Om^+$ be neither rotational nor a ruled
surface (see e.g. [Is1]). 
To reduce this problem to the boundary variational
inequality we need to reduce equivalently the nonhomogeneous
equation \refn{5.1} to the homogeneous one. To this end, consider the following
auxiliary problem:\\
Find a weak solution $U_0=(u_0,\om_0)^\top\!\in \![H^1(\Omega^+)]^6$
of the equation \refn{5.1} in $\Om^+$ satisfying on $S$ the
following conditions:
$$
\{u_{0n}\}^+=0,\quad \{\tau^{(n)}_s(U_0)\}^+=0,\quad\{\Mc
U_0\}^+=0\quad\text{on} \quad S.
$$
As it is known (see [GGN1], Theorem\,4.4) this problem is uniquely
solvable, since $S$ is neither rotational nor ruled. If $W\in [H^1(\Omega^+)]^6$ is a solution of
Problem $(B)$  and $U_0\in [H^1(\Omega^+)]^6$ is a solution of the above
auxiliary problem, then the difference $U:=W-U_0$ will be a
solution of the following problem.\\
{\it Problem\,$(B_0)$.} \textsl{Find a vector-function $U=(u,\om)^\top\in [H^1(\Omega^+)]^6$ which
is a weak solution of the homogeneous equation
\begin{equation}
 \label{5.2}
 L(\pa)\,U=0  \;\text{in} \;\Om^+,
\end{equation}
satisfying on $S$ the following conditions $\{\tau^{(n)}_s (U)\}^+\in [L_\infty (S)]^3$ and
\begin{eqnarray}
 \mbox{\rm (i)}
 &&
  \{\tau^{(n)}_n(U)\}^+=\psi ; \nonumber\\
  \mbox{\rm (ii)}
  &&
 \{\Mc U\}^+=\vfi; \nonumber \\
\mbox{\rm (iii)}
&&
\hskip-4mm
 \mbox{\rm(a)}
 \;\;
  \mbox{if} \;\;|\{\tau^{(n)}_s(U)\}^+|<g,\;\; \mbox{then}\;\; \{u_s\}^+=\vfi_0,\nonumber\\
&&
\hskip-4mm
 \mbox{\rm (b)}\;\;\mbox{if} \;\;|\{\tau^{(n)}_s(U)\}^+|=g, \;\; \mbox{then  there  exist nonnegative}\nonumber\\
 &&
  \hskip4mm
 \mbox{functions}\; \lambda_1\;  \mbox{and}\; \lambda_2\;
 \mbox{which do not vanish simultaneously and} \nonumber\\
 &&\hskip4mm \lambda_1\,\{u_s
\}^+=-\lambda_2\,\{\tau^{(n)}_s(U)\}^++\lambda_1\vfi_0,\nonumber
 \end{eqnarray}
where $\psi=F_0-\{\tau^{(n)}_n(U_0)\}^+$ and $\vfi_0=-\{u_{0s}\}^+.$}

Analogously to the previous coercive case (see Theorem\,3.2) it
can be shown that Problem\,$(B_0)$ is equivalent to the following boundary
variational inequality:\\
\textsl{Find $h_0=(h^{(1)}_0,h^{(2)}_0)^\top\in [H^{1/2}(S)]^6$ such that
the inequality
\begin{eqnarray}
\label{5.3}
 \langle\Ac^+h_0,h-h_0\rangle_{_{S}}+ j(h^{(1)})-j(h^{(1)}_0)\geq \langle
\psi,\,h^{(1)}_n-h^{(1)}_{0n} \rangle_{_S}+\langle\varphi,h^{(2)}-h^{(2)}_0\rangle_{_S}
\end{eqnarray}
holds for all $h=(h^{(1)},h^{(2)})^\top\in [H^{1/2}(S)]^6$, where now
$$
j(h^{(1)})=\int\limits_{S} g\,|h^{(1)}_s-\varphi_0|\,dS.
$$
}

Namely, the variational inequality \refn{5.3} and Problem\,$(B_0)$ are equivalent in the
following sense:
 If $U\in[H^{1}(\Omega^{+})]^6$ is a solution of the
Problem $(B_0)$, then $h=\{U\}^+\in [H^{1/2}(S)]^6$ is a solution
of the variational inequality \refn{5.3} and, vice versa, if
$h\in [H^{1/2}(S)]^6$ is a solution of the variational inequality
\refn{5.3}, then $G^+h\in[H^{1}(\Omega^{+})]^6$ is a weak
solution of Problem $(B_0)$. Here the operator $G^+$ is
defined by the equality \refn{3.17}.

Unfortunately, the variational inequality \refn{5.3} is not unconditionally solvable.

Now we derive the necessary condition of solvability of the
variational inequality \refn{5.3}.
Let $h_0=(h^{(1)}_0,h^{(2)}_0)^\top\in [H^{1/2}(S)]^6$ be a
solution of the inequality \refn{5.3}. Take
$h=(\vfi_0,h^{(2)}_0)^\top\in [H^{1/2}(S)]^6$ and
$h=(2h^{(1)}_0-\vfi_0,h^{(2)}_0)^\top$ instead of
$h=(h^{(1)},h^{(2)})^\top$
 in \refn{5.3} and take into account that the normal
component of $\vfi_0$ vanishes. We obtain
\begin{eqnarray}
\label{5.4}
 \langle\{\tau^{(n)}(G^+h_0)\}^+,h^{(1)}_0-\vfi_0\rangle_{_{S}}+ j(h^{(1)}_0)=
\langle\psi,\, h^{(1)}_{0n}\rangle{_{_S}}.
\end{eqnarray}
If we sum up the inequalities \refn{5.3} and \refn{5.4} we get
\begin{eqnarray}
\label{5.5}
&&
 \langle\{\tau^{(n)}(G^+h_0)\}^+,h^{(1)}-\vfi_0\rangle_{_{S}}+
 \langle\{\Mc(G^+h_0)\}^+,h^{(2)}-h^{(2)}_0\rangle_{_{S}}+j(h^{(1)})\nonumber\\
&&
\geq \langle \psi,\, h^{(1)}_{n}\rangle{_{_S}}+\langle\vfi,h^{(2)}-h^{(2)}_0\rangle_{_{S}}.
\end{eqnarray}
Rewrite \refn{5.5} as follows
\begin{eqnarray}
\label{5.6} &&
 \langle \psi,\, h^{(1)}_{n}\rangle{_{_S}}+\langle\vfi,h^{(2)}-h^{(2)}_0\rangle_{_{S}}
-\langle\{\tau^{(n)}(G^+h_0)\}^+,h^{(1)}-\vfi_0\rangle_{_{S}}\nonumber\\
&&-
 \langle\{\Mc(G^+h_0)\}^+,h^{(2)}-h^{(2)}_0\rangle_{_{S}}\leq
j(h^{(1)}).
\end{eqnarray}
Take here $2\vfi_0-h^{(1)}$ instead of $h^{(1)}$ and
$2h^{(2)}-h^{(2)}_0$ instead of $h^{(2)}$ to obtain
\begin{eqnarray}
\label{5.7} &&
 -\langle \psi,\, h^{(1)}_{n}\rangle{_{_S}}-\langle\vfi,h^{(2)}-h^{(2)}_0\rangle_{_{S}}
+\langle\{\tau^{(n)}(G^+h_0)\}^+,h^{(1)}-\vfi_0\rangle_{_{S}}\nonumber\\
&&+
 \langle\{\Mc(G^+h_0)\}^+,h^{(2)}-h^{(2)}_0\rangle_{_{S}}\leq
j(h^{(1)}).
\end{eqnarray}
From \refn{5.6} and \refn{5.7} we get the inequality
\begin{eqnarray}
\label{5.8}
 &&
 \Big|\langle\{\tau^{(n)}(G^+h_0)\}^+,h^{(1)}-\vfi_0\rangle_{_{S}}+
 \langle\{\Mc(G^+h_0)\}^+,h^{(2)}-h^{(2)}_0\rangle_{_{S}}\nonumber  \\
 &&
 - \langle \psi,\, h^{(1)}_{n}\rangle{_{_S}} -
\langle\vfi,h^{(2)}-h^{(2)}_0\rangle_{_{S}} \Big|
\leq j(h^{(1)})\quad \fa h=(h^{(1)},h^{(2)})^\top\in
[H^{1/2}(S)]^6.
\end{eqnarray}
Let $h^{(1)}-\vfi_0=\vartheta$ and $h^{(2)}-h^{(2)}_0=a$, where
$\vartheta=[a\times x]+b$ with arbitrary constant
vectors $a$ and $b$. Set $\chi:=(\vartheta,a)^\top$.
 Since $\langle\{T(\pa,n)(G^+h_0)\}^+,\chi\rangle_{_{S}}=0$, we have
from \refn{5.8}
\begin{eqnarray}
\label{5.9}
\Big| \langle \psi,\,\vte_{n}\rangle{_{_S}} +
\langle\vfi,a\rangle_{_{S}}\Big|\leq \int\limits_{S}
g\,|\vte_{s}|\,dS\quad\fa\chi=(\vartheta,a)^\top\in\Lambda(S),
\end{eqnarray}
where $\Lambda(S)=\ker \Ac^+$.  Thus, if
$h_0=(h^{(1)}_0,h^{(2)}_0)^\top\in [H^{1/2}(S_2)]^6$ is a solution
of the variational inequality \refn{5.3}, then \refn{5.9} holds
for all $\chi\in\Lambda(S)$, i.e., \refn{5.9} is a necessary
condition for solvability of \refn{5.3}.

Let us show that if \refn{5.9} holds with strict inequality sign,
i.e.,
\begin{eqnarray}
\label{5.10}
 \int\limits_{S}
g\,|\vte_{s}|\,dS-\Big|\langle \psi,\,\vte_{n}\rangle{_{_S}} +
\langle\vfi,a\rangle_{_{S}}\Big|>0
\quad\fa\chi=(\vartheta,a)^\top\in\Lambda(S),\;\,\chi\neq 0,
\end{eqnarray}
then the variational inequality \refn{5.3} is solvable. Since
$\Lambda(S)$ is a finite dimensional space,
\refn{5.10} can be sharpened to
\begin{eqnarray}
\label{5.11}
 \int\limits_{S}
g\,|\vte_{s}|\,dS-\Big|\langle \psi,\,\vte_{n}\rangle{_{_S}} +
\langle\vfi,a\rangle_{_{S}}\Big| \geq M \|\chi\|_{{[L_2(S)]^6}}
\quad \forall \chi=(\vartheta,a)^\top\in\Lambda(S),
\end{eqnarray}
with some positive constant $M>0$.

Therefore it suffices to show that \refn{5.10} is a sufficient condition
of solvability of the variational inequality \refn{5.3}. We proceed as follows.

Let $P$ be the operator of orthogonal projection (in the sense of
the space $[L_2(S)]^6$) of the space $[H^{1/2}(S)]^6$ on the space
$\Lambda(S)$ and $Q=I-P$. For any $h\in [H^{1/2}(S)]^6$ we have the representation
 $h=Qh+Ph$, where
$Ph=(\vte,a)^\top\in\Lambda(S),\;Qh=(\psi^{(1)},\psi^{(2)})^\top\in
[\Lambda(S)]^\bot$.

Consider the functional $\Jc$ on the space $[H^{1/2}(S)]^6$:
$$
\Jc
(h)=\frac{1}{2}\la\Ac^+h\,,\,h\ra_{_{S}}+j(h^{(1)})-
\langle\psi,\, h^{(1)}_{n}\rangle{_{_S}} - \langle\vfi\,,\,h^{(2)}\rangle_{_{S}},\quad
h=(h^{(1)},h^{(2)})^\top\in [H^{1/2}(S)]^6.
$$
Since the operator $\Ac^+$ is self-adjoint, as in the previous
case, the solvability of the inequality \refn{5.3} is equivalent
to the minimization problem for the functional $\Jc $ on the
space $[H^{1/2}(S)]^6$. Now we show that the functional $\Jc$ is
coercive, i.e.,
$$
\Jc
(h)\to+\infty\;\,\text{as}\;\,\|h\|_{_{[H^{1/2}(S)]^6}}\to\infty.
$$
Since $\Ac^+$ is self-adjoint, we have
$\la \Ac^+(Qh+Ph),Qh+Ph\ra_{_{S}}=\la \Ac^+Qh,Qh\ra_{_{S}}$ and from
Lemma 3.1 (e) we derive with some positive constant $c$,
\begin{align*}
 \Jc (h) & =\Jc (Qh+Ph)\\
&= \frac{1}{2}\, \la
\Ac^+Qh\,,\,Qh\ra_{_{S}}+\int\limits_{S}g\,|\psi^{(1)}_s+\vte_s-\vfi_0|\,dS-
\langle\psi,\,\psi^{(1)}_n+\vartheta_n\rangle{_{_S}}
 -\la\vfi,\psi^{(2)}+a\ra_{_{S}} \\
&\geq
c\,\|Qh\|_{_{[H^{1/2}(S)]^6}}^2-\int\limits_{S}g\,|\psi^{(1)}_s-\vfi_0|\,dS-
\langle\psi,\,\psi^{(1)}_n\rangle{_{_S}}
 -\la\vfi\,,\,\psi^{(2)}\ra_{_{S}} \\
&\hskip5mm+  \int\limits_{S}g\,|\vte_s|\,dS-
\langle \psi,\,\vte_{n}\rangle{_{_S}} -
\langle\vfi,a\rangle_{_{S}} \,.
\end{align*}

From this inequality, taking into account \refn{5.11} we obtain
for $\chi:=Ph$
\begin{eqnarray}
\label{5.12}
 \Jc (h)\geq
c\,\|Qh\|_{_{[H^{1/2}(S)]^6}}^2-c_1\,\|Qh\|_{_{[H^{1/2}(S)]^6}}+M\,\|\chi\|_{_{[L_2(S)]^6}}-c_2,
\end{eqnarray}
where $c_1$ and $c_2$ are positive constants.

It is easy to see that for $h=Qh+Ph$ the norm $\|h\|_{_{[H^{1/2}(S)]^6}}$ is
equivalent to the norm
$\|Qh\|_{_{[H^{1/2}(S)]^6}}+\|Ph\|_{_{[L_2(S)]^6}}$.
Therefore, from \refn{5.12} we see that, if
$\|h\|_{_{[H^{1/2}(S)]^6}}\rightarrow\infty$, then $\Jc (h)\rightarrow+\infty$
which proves the coercivity of the
functional $\Jc$. Due to the general theory of variational inequalities
(see \cite{GLT1}, \cite{DuLi1}) we conclude that the convex
continuous functional $\Jc(h)$ has a minimum on $[H^{1/2}(S)]^6$
and the minimizing function is a solution of the boundary
variational inequality \refn{5.3}.

Let $h,\,h^*\in[H^{1/2}(S)]^6$ be two arbitrary solutions of the
boundary variational inequality \refn{5.3}. It is easy to show
that then
$$
\la\Ac^+(h-h^*),h-h^*\ra_{_{S}}=0
$$
and consequently
$$
h-h^*=([a\times x]+b,a)\in \Lambda(S).
$$
Thus, from the above results we can formulate the
following assertion.
\begin{theorem}
Suppose, $S_1=\varnothing$. Let $X\in[L_2(\Omega^+)]^6,$ $F_0\in
L_\infty(S)$,  $\vfi\in [H^{-{1/2}}(S)]^3$ and  $g=\Fc|F_0|$ with
$\Fc : S\rightarrow[0,\infty)$ being a bounded measurable function,
   and  let \refn{5.11} hold. Then the boundary variational
inequality \refn{5.3} is solvable and the solutions are
determined modulo a generalized rigid displacement vector.
\end{theorem}
Due to the equivalence of the boundary variational inequality
\refn{5.3} and Problem (B) the counterpart of Theorem 5.1 holds
for Problem (B) as well.

Analogously to the semicoercive case, in the same way
we can investigate the contact Problem (C)
 when instead of the Dirichlet condition
\refn{3.3} the Neumann condition
$$
r_{_{S_1}}\{T(\pa,n) U\}^+=\Psi
$$
is given on the part $S_1$ of the boundary, where $\Psi\in
[\wt{H}^{-{1/2}}(S_1)]^6$; in contrast to the previous case,
here the vector $\vfi$ in the condition \refn{3.4} must be in
$[\wt{H}^{-{1/2}}(S_2)]^3$ and all other conditions on $S_2$
remain the same. Now, we need a solution to the following auxiliary problem:
Find a vector-function $U_0=(u_0,\om_0)^\top\in [H^1(\Omega^+)]^6$
which is a weak solution in $\Om^+$ of the equation
$$
 L(\pa)\,U_0+X=0
$$
and satisfies the following conditions:
$$
r_{_{S_1}}\{T(\pa,n) U_0\}^+=0\quad\text{on}\quad S_1
$$
$$
r_{_{S_2}}\{u_{0n}\}^+=0,\;\;r_{_{S_2}}\{\tau^{(n)}_s(U_0)\}^+=0,\;\;r_{_{S_2}}\{\Mc
U_0\}^+=0\;\;\text{on}\;\; S_2
$$
As it is known (see \cite{GGN1}, Theorem 4.4) if the part
$S_2$ is neither rotational nor ruled surface, this problem has a
unique solution.

Let $W\in [H^1(\Omega^+)]^6$ be a solution of the above Problem (C)
and $U_0\in [H^1(\Omega^+)]^6$ be a solution of the auxiliary
problem, then the difference $U:=W-U_0$ will be a solution of the
following problem.\\
{\it Problem} $(C_0)$. \textsl{Find a weak solution $U=(u,\om)^\top\in
[H^1(\Omega^+)]^6$ of the equation
$$
 L(\pa)\,U=0\quad \text{in}\quad \Omega^+,
$$
satisfying the inclusion $r_{_{S_2}}\{\tau^{(n)}_s(U)\}^+\in
[L_\infty (S_2)]^3$ and the following conditions:
\begin{eqnarray*}
 \mbox{\rm (i)}
 &&
   r_{_{S_1}}\{T(\pa,n) U\}^+=\Psi\quad  \mbox{on}\quad S_1;\\
 \mbox{\rm (ii)}
  &&
   r_{_{S_2}}\{\tau^{(n)}_n(U)\}^+=\psi\quad
 \mbox{on} \quad S_2;\\
 \mbox{\rm (iii)}
 &&
   r_{_{S_2}}\{\Mc U\}^+=\vfi\quad  \mbox{on}\quad S_2;\\
\mbox{\rm (iv)}
&&
\hskip-4mm
 \mbox{\rm (a)} \;\;
\mbox{if} \;\;
 |r_{_{S_2}}\{\tau^{(n)}_s(U)\}^+|<g,\;\;
\mbox{then} \;\; r_{_{S_2}}\{u_s\}^+=\varphi_0\quad
 \mbox{on}\;\,
S_2,\\
&&
\hskip-4mm
\mbox{\rm (b)}\; \;\mbox{if} \;\;
|r_{_{S_2}}\{\tau^{(n)}_s(U)\}^+|=g, \;\; \mbox{then  there  exist nonnegative}\nonumber\\
 &&
 \hskip3mm
  \mbox{ functions}\; \lambda_1\;  \mbox{  and}\; \lambda_2\;
 \mbox{which do not vanish simultaneously and} \nonumber\\
 &&
 \hskip3mm
 \; \lambda_1\,r_{_{S_2}}\{u_s
\}^+=-\lambda_2\,r_{_{S_2}}\{\tau^{(n)}_s(U)\}^++\lambda_1\varphi_0
\quad  \mbox{on} \;\,S_2,
 \end{eqnarray*}
 where $g$ is defined  by  formula  \refn{3.2},
 $\psi=F_0-r_{_{S_2}}\{\tau^{(n)}_n(U_0)\}^+$ and  $\varphi_0=-r_{_{S_2}}\{u_{0 s}\}^+$.
}\\
The equivalent boundary variational inequality to this problem is
the following:
\textsl{Find $h_0=(h^{(1)}_0,h^{(2)}_0)^\top\in [H^{1/2}(S)]^6$ such that
the inequality
\begin{eqnarray}
\label{5.13}
&&
 \langle\Ac^+h_0,h-h_0\rangle_{S}+ j(h^{(1)})-j(h^{(1)}_0)\geq\langle\psi,
 r_{_{S_2}}(h^{(1)}_n-h^{(1)}_{0n})\rangle_{S_2} \nonumber\\
 &&
 +\langle r_{S_2}\varphi,r_{S_2}(h^{(2)}-h^{(2)}_0)\rangle_{S_2}
 +\langle r_{S_1}\Psi,r_{S_1}(h-h_0)\rangle_{S_1}
\end{eqnarray}
holds for all $h=(h^{(1)},h^{(2)})^\top\in [H^{1/2}(S)]^6$.
}

The variational inequality \refn{5.13} and Problem (C$_0$)
are equivalent in the sense described in Theorem 3.2.
The necessary condition of solvability of the variational inequality \refn{5.13}
reads
\begin{eqnarray}
\label{5.14}
 \left|\langle\psi,r_{_{S_2}}\{\vte_{n}\}^+\rangle_{_{S_2}}+
\langle r_{S_2}\vfi,a\rangle_{_{S_2}}+\langle
r_{S_1}\Psi,r_{_{S_1}}\{\chi\}^+\rangle_{_{S_1}}\right |\leq
\int\limits_{S_2} g\,|\{\vte_{s}\}^+|\,dS,
\end{eqnarray}
where
$\chi=(\vartheta,a)^\top\in\Lambda(S),\;\,\vartheta=[a\times
x]+b,\;\,a,\,b\in\R^3.$

If \refn{5.14} holds with the strict inequality sign, then
we can sharpen it to
\begin{eqnarray}
\label{5.15}
\hskip-5mm
&
\dst
 \int\limits_{S_2}\!
g\,|\{\vte_{s}\}^+|dS\!-\!\left|\langle\psi,r_{_{S_2}}\{\vte_{n}\}^+\rangle_{_{S_2}}\!\!\!+\!
\langle r_{_{S_2}}\vfi,a\rangle_{_{S_2}}\!\!\!+\!\langle
r_{_{S_1}}\Psi,r_{_{S_1}}\{\chi\}^+\rangle_{_{S_1}}\right|
\! \geq \! M\|\chi\|_{[L_2(S)]^6} \\
\hskip-5mm
&
\fa\chi=(\vartheta,a)^\top\in\Lambda(S)\nonumber
\end{eqnarray}
with some positive constant $M$, since $\Lambda(S)$ is a finite dimensional space.
By the same arguments as above it can be shown that this condition is
sufficient for the solvability of the
variational inequality \refn{5.13}. Finally we arrive at the following theorem.
\begin{theorem}
If $\Psi\in [\wt{H}^{-{1/2}}(S_1)]^6,$ $\vfi\in
[\wt{H}^{-{1/2}}(S_2)]^3,$ $\psi\in \wt{H}^{-{1/2}}(S_2),$ $\Fc:S_2\rightarrow[0,\infty)$ is a
bounded measurable function, $g=\Fc|F_0|$ and \refn{5.15} holds,
then there exists a solution $h_0$ of the variational inequality
\refn{5.13} and $G^+h_0$  solves Problem\,$(C_0)$. Solutions of the variational inequality \refn{5.13} and Problem\,$(C_0)$ are defined modulo generalized rigid
displacement vectors.
\end{theorem}

\section{Exterior problems}

First of all let us observe that the bilinear form
$$
B(U,V):=\int\limits_{\Om^-}E(U,V)\,dx
$$
is well defined for vectors
$U=(u,\om)^\top\in[H^{1}_{loc}(\Omega^{-})]^6$ and
$V=(v,w)^\top\in[H^{1}_{loc}(\Omega^{-})]^6$ satisfying the decay
conditions (Z) at infinity.

Let $X\in[L_{2,comp}(\Omega^-)]^6,$
$f\in [H^{{1/2}}(S_1)]^6,$ $F_0\in L_\infty(S_2)$, $\vfi\in
[H^{-{1/2}}(S_2)]^3$, and $g=\Fc|F_0|$ with $\Fc:\;S_2\rightarrow[0,\infty)$
 being a bounded measurable function.
Consider the following bilateral contact problem with friction.\\
{\it Problem}  $(D)$. \textsl{ Find a weak solution $U=(u,\om)^\top\in
[H_{\loc}^1(\Omega^-)]^6$ of equation
\begin{eqnarray}
 \label{6.1}
 L(\pa)\,U+X=0\quad \text{in}\quad \Omega^-,
\end{eqnarray}
satisfying the decay conditions $(Z)$ at infinity, the inclusion
$r_{_{S_2}}\{\tau^{(n)}_s(U)\}^-\in [L_\infty (S_2)]^3$ and the boundary conditions
\begin{eqnarray}
\label{6.2} \mbox{\rm (i)}
&&
  r_{_{S_1}}\{U\}^-=f\;\;  \mbox{on}\;\; S_1;\\
 \mbox{\rm (ii)}
 &&
  r_{_{S_2}}\{\tau^{(n)}_n(U)\}^-=F_0\;\; \mbox{on} \;\; S_2; \nonumber\\
\label{6.3} \mbox{\rm (iii)}
&&
 r_{_{S_2}}\{\Mc U\}^-=\vfi\;\; \mbox{on}\;\;S_2;\\
 \mbox{\rm (iv)}
 &&
 \hskip-4mm
 \mbox{\rm (a)}\;\;
\mbox{if} \;\;|r_{_{S_2}}\{\tau^{(n)}_s(U)\}^-|<g,\;\;
\mbox{then}\;\; r_{_{S_2}}\{u_s\}^-=0,\nonumber\\
&&
 \hskip-4mm
 \mbox{\rm (b)}\;\;
\mbox{\rm if} \;\;|r_{_{S_2}}\{\tau^{(n)}_s(U)\}^-|=g, \;\; \mbox{then  there  exist nonnegative}\nonumber\\
 &&\hskip4mm
 \mbox{functions}\; \lambda_1\;  \mbox{  and}\; \lambda_2\;
 \mbox{which do not vanish simultaneously and} \nonumber\\
 &&
 \hskip4mm \lambda_1\,r_{_{S_2}}\{u_s
\}^-=-\lambda_2\,r_{_{S_2}}\{\tau^{(n)}_s(U)\}^-\;\;\mbox{on}\;\; S_2.  \nonumber
 \end{eqnarray}}
To reduce this problem to the boundary variational inequality, as a first step, again we
have to reduce the nonhomogeneous equation \refn{6.1}
and the nonhomogeneous Dirichlet condition \refn{6.2} to the homogeneous ones. For this
purpose consider the following auxiliary problem:
Find a weak solution $U_0\in [H_{\loc}^1(\Omega^-)]^6$ of the equation
$$
 L(\pa)\,U_0+X=0\quad \text{in}\quad \Omega^-
$$
satisfying the decay conditions (Z) at infinity and the boundary
conditions on $S$:
\begin{eqnarray*}
    r_{_{S_1}}\{U_0\}^-=f\quad  \mbox{\rm on}\quad S_1,\qquad
    r_{_{S_2}}\{T(\pa,n)U_0)\}^-=0\quad
 \mbox{\rm on} \quad S_2.
\end{eqnarray*}
By the same approach as in the case of the interior problem, with the help of the solution vector
$U_0=(u_0, \omega_0)^\top$ the original Problem $(D)$ can be reduced
to the following one.\\
{\it Problem}  $(D_0)$. \textsl{Find a weak solution $U=(u,\om)^\top\in
[H_{\loc}^1(\Omega^-)]^6$ of the homogeneous equation
\begin{eqnarray}
 \label{6.4}
   L(\pa)\,U=0\quad \text{in}\quad \Omega^-
  \end{eqnarray}
satisfying the decay conditions $(Z)$ at infinity, the inclusion
 $r_{_{S_2}}\{\tau^{(n)}_s(U)\}^-\in [L_\infty (S_2)]^3$
and boundary conditions
\begin{eqnarray*}
\mbox{\rm (i)}
&&
  r_{_{S_1}}\{U\}^-=0\;\; \mbox{on}\;\; S_1;\\
\mbox{\rm (ii)}
&&
  r_{_{S_2}}\{\tau^{(n)}_n(U)\}^-=F_0\;\; \mbox{on} \;\;  S_2;\\
\mbox{\rm (iii)}
&&
 r_{_{S_2}}\{\Mc
U\}^-=\vfi\;\; \mbox{on}\;\;S_2;\\
 \mbox{\rm (iv)}
 &&
 \hskip-4mm
 \mbox{\rm (a)}\;\;
\mbox{if} \;\;|r_{_{S_2}}\{\tau^{(n)}_s(U)\}^-|<g,\;\;
\mbox{then}\;\; r_{_{S_2}}\{u_s\}^-=\varphi_0,\\
&&
\hskip-4mm
\mbox{\rm (b)}\;\;
\mbox{ if} \;\;|r_{_{S_2}}\{\tau^{(n)}_s(U)\}^-|=g, \;\; \mbox{then there  exist nonnegative}\nonumber\\
 &&
 \hskip4mm
 \mbox{functions}\; \lambda_1\;  \mbox{and}\; \lambda_2\;
 \mbox{which do not vanish simultaneously and} \nonumber\\
 &&
 \hskip4mm \lambda_1\,r_{_{S_2}}\{u_s \}^-=
 -\lambda_2\,r_{_{S_2}}\{\tau^{(n)}_s(U)\}^- + \lambda_1\,\varphi_0\;\; \mbox{on}\;\; S_2,
 \end{eqnarray*}
where
 $\varphi_0=-r_{S_2}\{u_{0s}\}^+$.}

To reduce this exterior problem equivalently to the boundary
variational inequality we apply the following representation of
solution of the equation \refn{6.4} from the class $[H_{\loc}^1(\Omega^-)]^6$ satisfying the
decay condition (Z) at infinity (see \cite{NGS1})
$$
U(x)=(G^-h)(x):=V(\Hc^{-1}h)(x)=\int\limits_S\Gamma(x-y)\,(\Hc^{-1}\,h)(y)\,dS_y,\quad
x\in\Om^-,
$$
where $\Gamma$ is the fundamental solution of the operator
$L(\pa)$, $\Hc$ is defined by \refn{3.11} and $h=\{U\}^-\in [H^{{1/2}}(S)]^6$.

Define the Steklov-Poincar\'{e} type operator
$\Ac^-$ with the help of the formula
$$
\Ac^-h:=\{T(\pa,n)(G^-h)\}^-\equiv\{T(\pa,n)V(\Hc^{-1}h)\}^-.
$$
Due to the properties of the single layer potential this operator
can be represented as
$$
\Ac^-=\left(2^{-1}I_6+\Kc\right)\Hc^{-1},
$$
where $\Kc$ is defined by \refn{3.16}.

The operator $\Ac^-$ possesses almost the same properties as $\Ac^+$. Namely,
\begin{eqnarray*}
&&
 \mbox{\rm (a)}\;\;\langle\Ac^-h',h''\rangle_{S}=\langle\Ac^-h'',h'\rangle_{S} \quad
 \fa h'\in[H^{1/2}(S)]^6\;\text{and}\;\fa h''\in[H^{1/2}(S)]^6;\\
&& \mbox{\rm (b)}\;\;\Ac^-:[H^{1/2}(S)]^6\rightarrow
 [H^{-1/2}(S)]^6 \text{is a continuous operator};\\
 &&
 \mbox {\rm (c)} \;\text{there is a  constant} \,c_0>0,
 \text{such that}
 \langle\Ac^-h,h\rangle_{S}\geq
 c_0\|h\|^2_{[H^{1/2}(S)]^6}\;\fa
 h\in[H^{1/2}(S)]^6.
\end{eqnarray*}
Further, let us recall that
$
{\mathbb{H}}=\{h=(h^{(1)},h^{(2)})^\top\in
[H^{1/2}(S)]^6\,:\,r_{S_1}h=0\}
$
and consider the variational inequality on ${\mathbb{H}}$:
  \textsl{Find a vector function $h_0=(h^{(1)}_0,h^{(2)}_0)^\top\in {\mathbb{H}} $, such
 that the inequality
\begin{eqnarray}
\label{6.5}
 \langle\Ac^-h_0,h-h_0\rangle_{S}+ j(h^{(1)})-j(h^{(1)}_0)\geq\int\limits_{S_2}
F_0\,(h^{(1)}_n-h^{(1)}_{0n})\,dS+\langle\varphi,r_{S_2}(h^{(2)}-h^{(2)}_0)\rangle_{S_2},
\end{eqnarray}
holds for all $h=(h^{(1)},h^{(2)})^\top\in {\mathbb{H}} $.
 Here $j(\cdot)$ is defined by the relation \refn{3.20}.}

Using the same arguments as in Theorems 3.2 and 4.1 we can prove
the equivalence of the variational inequality \refn{6.5} and  Problem\,$(D_0)$, and
the uniqueness theorem of solution to the variational inequality \refn{6.5}. To prove the
existence of solutions, we consider the following functional on
the closed subspace $\mathbb{H}$
\begin{eqnarray*}
 \Jc(h)=\frac{1}{2}\,\langle\Ac^-h,h\rangle_{S}+ j(h^{(1)})-\int\limits_{S_2}
F_0\,h^{(1)}_n\,dS-\langle\varphi,r_{S_2}h^{(2)}\rangle_{S_2},\;\;\;
  h=(h^{(1)},h^{(2)})^\top\in \mathbb{H}.
\end{eqnarray*}
Due to the symmetry property (a) of the operator $\Ac^-$, the
solvability of the variational inequality \refn{6.5} is equivalent to
the minimization problem for the functional $\Jc (h)$ on the set
${\mathbb{H}}$ . It is easy to see that from the coercivity property (c) of
the operator $\Ac^-$ and from the inequality $j(h^{(1)})\geq 0$ we
obtain the estimate from below for the functional $\Jc(h)$
$$
\Jc (h)\geq
c_1\|h\|^2_{[H^{1/2}(S)]^6}-c_2\|h\|_{[H^{1/2}(S)]^6}\quad \fa
h\in {\mathbb{H}} .
$$
Hence the coercivity of the functional $\Jc(h)$ follows,
i.e., $\Jc (h)\rightarrow+\infty$ if
$\|h\|_{_{[H^{1/2}(S)]^6}}\rightarrow\infty$. Due to the theory of
variational inequalities (see \cite{GLT1}, \cite{DuLi1}) we
conclude that the convex continuous functional $\Jc(h)$ has a
minimum on $[H^{1/2}(S)]^6$ and the minimizing function  $h_0$ is a
unique solution of the variational inequality \refn{6.5}. Consequently, the unique solution of Problem (D$_0$)
can be represented in the form $U=G^-h_0$. Finally, we arrive at the following existence result.
\begin{theorem}
Let $X\in[L_{2,\,comp}(\Omega^-)]^6,$ $\vfi\in [H^{-{1/2}}(S_2)]^3,$ $F_0\in
L_\infty(S_2),$ $f\in[H^{1/2}(S_1)]^6$ and
$\Fc:\;S_2\rightarrow[0,\infty)$ be a bounded measurable function.
Then Problem $(D)$ has a unique solution in the space
$[H^1_{\loc}(\Om^-)]^6$ satisfying the decay conditions $(Z)$ at
infinity.
\end{theorem}
\begin{remark}
By the same approach one can investigate the problem when either

\item $(a)$  $S_1=\varnothing$ and the friction
conditions are considered on the whole boundary\\
or
 \item $(b)$ $S_1\neq\varnothing$ and instead of the Dirichlet condition \refn{6.2}
 there is given the Neumann condition
$$
r_{_{S_1}}\{T(\pa,n) U\}^+=\Psi,
$$
 where $\Psi\in[\wt{H}^{-{1/2}}(S_1)]^6$. Note that the vector $\vfi$ involved
 in the boundary condition \refn{6.3} now should be from the space $[\wt{H}^{-{1/2}}(S_2)]^3$.\\
 Both these problems are uniquely solvable.
\end{remark}

\end{document}